\documentclass[11pt,reqno]{amsart}

\usepackage[english]{babel}
\usepackage[utf8x]{inputenc}
\usepackage{amsmath}
\usepackage{amsthm}
\usepackage{graphicx}
\usepackage{float}
\usepackage[colorinlistoftodos]{todonotes}
\usepackage{amsmath,amssymb,mathrsfs,amsthm,dsfont}
\usepackage[T1]{fontenc}
\usepackage{tikz}
\usetikzlibrary{matrix}
\title[Simpson index for Moran process in random environment]{On the Simpson index for the Moran process with random selection and immigration }
\author[A. Guillin]{\textbf{\quad {Arnaud} Guillin $^{\diamondsuit}$ \, \, }}
\address{{\bf {Arnaud} GUILLIN}\\ Laboratoire de Math\'ematiques Blaise Pascal, CNRS UMR 6620, Universit\'e Clermont-Auvergne,
avenue des Landais, F-63177 Aubi\`ere.} \email{arnaud.guillin@uca.fr}
\author[F. Jabot]{\textbf{\quad {Franck} Jabot  $^{\clubsuit}$ \, \, }}
\address{\textbf{{Franck} JABOT}\\ Laboratoire d'Ing\'eni\'erie pour les Syst\`emes Complexes, IRSTEA, Campus des C\'ezeaux
9, avenue Blaise Pascal - CS 20085
63178 Aubi\`ere}
\email{franck.jabot@irstea.fr}
\author[A. Personne]{\textbf{\quad {Arnaud} Personne $^{\diamondsuit}$ \, \, }}
\address{{\bf {Arnaud} PERSONNE}\\ Laboratoire de Math\'ematiques Blaise Pascal, CNRS UMR 6620, Universit\'e Clermont-Auvergne,
avenue des Landais, F-63177 Aubi\`ere.} \email{arnaud.personne@uca.fr}
\theoremstyle{definition}

\newtheorem{remarque}{Comment}

\theoremstyle{plain}
\newtheorem{theorem}{Theorem}
\newtheorem{lemma}[theorem]{Lemma}
\newtheorem{prop}[theorem]{Proposition}

\newcommand{\E}{\mathbb{E}}

\newcommand{\pp}{\mathbb{P}}

\begin{document}
\parindent=0pt

\maketitle

 \begin{center}

\textsc{$^{\diamondsuit}$ Universit\'e Clermont-Auvergne}
\smallskip

\textsc{$^{\clubsuit}$ Irstea}
\smallskip

\end{center}

\begin{abstract}
Moran or Wright-Fisher processes are probably the most well known model to study the evolution of a population under various effects. Our object of study will be the Simpson index which measures the level of diversity of the population, one of the key parameter for ecologists who study for example forest dynamics. Following ecological motivations, we will consider here the case where there are various species with fitness and immigration parameters being random processes (and thus time evolving). To measure biodiversity, ecologists generally use the Simpson index, who has no closed formula, except in the neutral (no selection) case via a backward approach, and which is difficult to evaluate even numerically when the population size is large. Our approach relies on the large population limit in the "weak" selection case, and thus to give a procedure which enable us to approximate, with controlled rate, the expectation of the Simpson index at fixed time. Our approach will be forward and valid for all time, which is the main difference with the historical approach of Kingman, or Krone-Neuhauser. We will also study the long time behaviour of the Wright-Fisher process in a simplified setting, allowing us to get a full picture for the approximation of the expectation of the Simpson index. 
\end{abstract}

\bigskip

\bigskip

 
\tableofcontents


\section{Introduction}


Community ecology has been deeply shaked by the book of Hubbell (2001) \cite{hubbel} that elaborated on the idea that the dynamics of ecological communities might be mainly shaped by random processes. A number of predictions made by this neutral theory of biodiversity have been indeed corroborated by empirical evidence (Hubbell \cite{hubbel}, Condit et al.\cite{Conditall}, Jabot and Chave \cite{JabotChave2009}). This good performance of neutral models for reproducing empirical patterns has stimulated the mathematical study of neutral ecological models (Etienne \cite{etienne2005}, Fuk et al \cite{FOC2017}), in connection with the rich body of work on evolutionary neutral models (Volkov et al. \cite{VBHM2003}, Ewens \cite{ewens2004}, Muirhead and Wakeley \cite{mw2009}).\\ 
More recent empirical evaluations of neutral predictions on tropical forest data have focused on the temporal dynamics of individual populations and have shown that the temporal variance of population sizes was actually larger than the one typically predicted by neutral models (Chisholm et al. \cite{Chisholmall}). These authors have suggested that this may be due to species-specific responses to the temporal variability of the environment. Subsequent modelling studies have elaborated on this idea (Kalyuzhny et al. \cite{kalyuzhny}, Jabot and Lohier \cite{jabotlohier}, Krone and Neuhauser \cite{KNeuhauser} and more mathematically-oriented contributions on this topic have since been made (\cite{GriffthsTavare,expansionmultiallelicWF,daninoshnerb,danino2016effect,danino2018stability,danino2018fixation,fung2017species,Viotprocees,Grieshammer,cmv}). All these studies share the same idea that community dynamics is influenced by a species-specific selection coefficient and that this selection coefficient is temporally varying, so that good and bad periods are subsequently experienced by all species within the community.\\
The aim of this contribution is to provide a comprehensive study of a simple model encapsulating this kind of ideas within a weak selection framework. Subsequent works within a strong selection framework will complement the present study.

Indeed, the simplest model for the evolution of population is surely the Moran model (or its cousin the discrete Wright-Fisher model), in which in a given population an individual is chosen to die (uniformly in the population) and then a child chooses his parent proportionally to the abundance in the previous population. One may also add immigration, i.e. a probability that the child comes from another community, and selection so that some species (or traits) have a selective advantage. Kalyuzhny et al \cite{kalyuzhny}, to bypass the neutrality of the model, chosed to consider immigration and selection as random processes (independent of the Moran system), but which preserves neutrality "in mean". One of the main goal is of course to study the effect of these models on biodiversity and its evolution. There are many ways to measure biodiversity. We will focus here on the Simpson index \cite{Simpson}, usually considered in neutral model \cite{etienneOlff2}: it measures the  probability that two individuals uniformly chosen may be of the same species. More precisely denoting $X^{i}$ the number of individuals of species $i$, $S$ the number of species and $J$ the toal size of the population, the Simpson index is given by
$${\mathcal S}=\sum\limits_{i=1}^{S}\frac{X^{i}(X^{i}-1)}{J(J-1)}$$
thus varying (roughly) from $0$ to $1$, from maximal to minimal diversity. Using backward approach, and Kingman's coalescent, an explicit formula may be given for the (asymptotic in time) Simpson index in the neutral case with immigration as $\frac{1-m}{1+J-2m}$, as otherwise it is 1 as a particular specie will almost surely invade all the population. Note that a closed formula (even for the expectation) of the Simpson index at a given time, is usually not reachable.\\
We will consider in this paper this Moran model with immigration and selection as general random process. As said previously such models were recently considered for example by Griffiths \cite{Griffiths,GriffthsTavare}, Kalyuznhy et al \cite{kalyuzhny} for a simulation study, but no theoretical framework towards the Simpson index. The backward approach constitutes the works of Krone and Neuhauser \cite{KNeuhauser,NKrone} leading to new coalescent type processes which are however quite difficult to study and may not give a closed formula for the Simpson index. A recent work by Grieshammer \cite{Grieshammer} considers a forward in time approach but he does not focus on the Simpson index. Our approach is only forward here. As is often done in population genetics, we will consider the large population approximation. Our first task is then to justify this asymptotic to a Wright-Fisher diffusion process in random environment in the weak immigration and selection case. As a flavour, with only two species,  the evolution of the proportion of one species is given by
$$dX_t=m_t(p_t-X_t)dt+s_tX_t(1-X_t)dt+\sqrt{2X_t(1-X_t)}dB_t$$
where $m_t$ is the immigration process, $p_t$ the probability that this species is chosen, and $s_t$ the selection advantage. It will be done by the usual martingale method. Another quantitative approach will be considered in \cite{ggp2018}. The Simpson index is then a quadratic form involving the proportion of each species, and by It\^o's formula it involves higher order term. The equation for the expectation of the Simpson index is therefore not closed. We will then introduce a quantitative approximation procedure for the expectation of the Simpson index, in the quenched case (corresponding to a given random environment) and in some particular case in the annealed case for two or more species. We will also study in some simple case (constant parameter) the long time behaviour of the Simpson index. It is reminiscent with the very recent work of Coron, M\'el\'eard and Villemonais \cite{cmv} (discovered while finishing this work). Very schematically our approach is the following
\begin{enumerate}
\item approximate the true discrete process by a SDE, i.e. Wright-Fisher process;
\item approximate the expectation of the Simpson index for the Wright-Fisher process by a deterministic ODE;
\item approximate the infinite time expectation of the Simpson index by the finite time, through evaluation of the speed of convergence towards equilibrium of the Wright-Fisher process.
\end{enumerate}
In Section 2, we introduce the Moran model in random environment and prove its convergence in the large population limit. In Section 3, we focus on the two species case where the approximation of the Simpson index is studied in the quenched case as well as its long time behaviour. Section 4 generalizes to a large number of species and also considers the annealed case when the selection parameter has a particular form (a variant of a Wright-Fisher process). The last section contains technical proofs or recall some known results for the Wright-Fisher process.

\section{The Moran model in random environment and its approximation in large population }

\subsection{Discrete model with selection and immigration}
\quad\\
In this section we describe in detail the discrete model, i.e. the Moran process, which is the basic of our study. One may also consider here the Wright-Fisher discrete process with adequat change.
The Moran process is an evolution of population model,  in which a single event occurs at each time step. More precisely each event corresponds to the death of an individual and the birth of another who replaces it.\\
We consider a population, whose size is constant over time equal to $J$, composed of $S+1$ species . The proportion of the $i$ species at the $n^{th}$ event is denoted $X^{i}_{n}$, $i \in \mathbb{S}=\{1,...,S+1\}$, $n\in \mathbb{N}$.\\As usual one we know $(X_n^{i})_{i=1,..,S}$, we deduce the proportion for the last species, $X_{n}^{S+1}=1-\sum\limits_{i=1}^{S}X^{i}_{n}$\\
We note $X_{n}$ the species vector or abundance vector  having for coordinate $i$, $X_{n}^{i}$.
The dynamics of evolution follows the following pattern at the step $n$:
\begin{enumerate}
\item The individual designated to die is chosen uniformly among the community.
\item The one which replaces it, chooses  his parent in the community with probability $1-m_{n}$ (filiation) or a parent from the immigration process with probability $m_{n}$ (immigration). The quantity $m_{n}$ varies between $0$ and $1$, it can be random and time dependent.
\item If there is immigration the chosen parent is from the species $i$,  with probability $p^{i}_{n}$ $i \in \mathbb{S}$. The $p_{n}^{i}$ verify $\sum\limits_{i \in \mathbb{S}}p^{i}_{n}=1$ and can be time dependent and random. We note $p$ for the vector having for coordinate $i$, $p_{n}^{i}$.
\item 
In a filiation, the chosen parent is of the species $i$ with probability $\frac{X^{i}_{n}(1+s^{i}_{n})}{1+\sum\limits_{k=1}^{S+1}X^{k}s^{k}_{n}}$. \\
The $s_{n}^{i}, i \in \mathbb{S}$ are the selection parameters , they may be time dependent and random.
 Furthermore, we  assume $s_{n}^{S+1}=0$. Indeed, we can obtain it from any configuration by changing all the coefficients by  $s_{n}^{i}=\frac{\tilde{s}^{i}_{n}-\tilde{s}^{S+1}_{n}}{1+\tilde{s}_{n}^{S+1}} $.
\end{enumerate}
We will assume throughout this work that $m_n,p_n,s_n$ are autonomous, in the sense that their evolution do not depend on $(X_n)_{n\ge0}$. We will further assume that $(m_n,p_n,s_n)_{n\ge0}$ is a Markov chain. Note also that $m_n,p_n,s_n$ may also depend in some sense of the size of the population $J$, but we do not add another superscript to get lighter notations.\\

This model therefore describes a Markovian dynamic in which selection and immigration play an important role. Immigration already introduced by Hubbell  \cite{hubbel} avoids the definitive invasion of the community by a species. Selection changes the dynamics of a species related to the neutral model( \cite{kalyuzhny}).
The temporal evolution of the population could be simulated numerically from the transition matrix of the Markov system. Let us describe precisely these transition probabilities for the evolution of proportions. The assumption for the dynamics for the immigration and selection will be given later on.\\
Let $x$ be the vector having for coordinate $i$, $x^{i}$  and suppose $m_{n}$ ,$p_{n}$  known.
Denote $\Delta=\frac{1}{J}$, so for the $i$ species:

    \begin{eqnarray*}
     P_{x^{i}+}&=& \pp(X^{i}_{n+1}=x+\Delta|X_{n}=x)\\ & =&(1-x^{i})\left(m_{n}p^{i}_{n}+(1-m_{n})\frac{x^{i}(1+s^{i}_{n})}{1+\sum\limits_{k=1}^{S+1}x^{k}s^{k}_{n}}\right),
     \end{eqnarray*}
     \begin{eqnarray*}
  P_{x^{i}-}   &=&\pp(X^{i}_{n+1}=x-\Delta|X_{n}=x))\\ & =&x^{i}\left(m_{n}(1-p^{i}_{n})+(1-m_{n})\left(1-\frac{x^{i}(1+s^{i}_{n})}{1+\sum\limits_{k=1}^{S+1}x^{k}s^{k}_{n}}\right)\right),
\end{eqnarray*}
     \begin{eqnarray*}
    P_{x^{i+}x^{j-}} &=&\mathbb{P}(\{X^{i}_{n+1}=x^{i}+\Delta\} \cap\{ X^{j}_{n+1}=x^{j}-\Delta\}|X_{n}=x)\\
     &=&x^{j}\left(m_{n}(1-p_{n}^{i})+(1-m_{n})\frac{x^{i}(1+s_{n}^{i})}{1+\sum\limits_{k=1}^{s+1}x^{k}s_{n}^{k}}\right).
\end{eqnarray*}

With the dynamics of $(m_n,p_n,s_n)$ given, one may of course simulate exactly the vector of proportion $X_n$ and thus evaluate the expectation of the Simpson index, which is what we will do to validate our approximation procedure, but when the population size $J$ is very large, it may be computationally too costly (and even impossible).
Thus we will approach the dynamics of this model by a stochastic differential system continuous in time.

\subsection{To a limit in large population}
\quad\\
In this section, we explain how approaching the dynamics of the preceding model  by a diffusion, and associated process for the immigration and selection processes, when $J$ goes to infinity. \\

We need to define a $J$ dependent time scale. Indeed, when $J$ goes to infinity, the time scale has to change, expectation and variance are about $\frac{1}{J}$ and $\frac{1}{J^2}$, and goes to  $0$ when $J $ goes to infinity. It corresponds to considering a large number of event for the Markov chain, to obtain a non-trivial convergence of our discrete process towards a limit process, i.e. not look at the event-by-event evolution as we did before but in packets of several events. \\
Several choices for scales are possible, each one leads to study a different process. We choose to study a continuous multidimensional diffusion in time.\\

\subsubsection{Diffusion approximation}
\quad\\
In a general framework, the limiting process we obtain is characterized by the first moment  and the covariance matrix of the infinitesimal variation of abundance.  \\
More precisely if we note  $\Delta_{t}$ the infinitesimal variation in time (which depends on the scale ) and $\Delta X_{t} =X_{t+\Delta_{t}}-X_{t}$ the infinitesimal variation of abundance, the diffusion process is characterized by the quantities:

$$b(x)=\lim\limits_{\Delta_{t}\Rightarrow 0}\frac{E[\Delta X_{t}|X_{t}=x]}{\Delta_{t}}$$
$$\sigma_{i,j}(x)=\lim\limits_{\Delta_{t}\Rightarrow 0}\frac{\mbox{Cov}[\Delta X_{t+\Delta_{t}}(i),\Delta X_{t+\Delta_{t}}(j)|X_{t}=x]}{\Delta_{t}} \quad i,j \in \mathbb{S}$$
\vspace{3mm}

The following property characterizes the order (relative to J) of the expectation, variance, and covariance of the abundance variation of a species during an event:

\begin{prop} 
\begin{enumerate}
\item $\E[X^{i}_{n+1}-X^{i}_{n}|X_{n}=x]=\Delta (P_{x^{i}+}-P_{x^{i}-})$
\item $\mbox{Var}[X^{i}_{n+1}-X^{i}_{n}|X_{n}=x]=\Delta^{2}((P_{x^{i}+}+P_{x^{i}-})-(P_{x^{i}+}-P_{x^{i}-})^{2})$
\end{enumerate}
\end{prop}
The proof is standard calculus and thus omitted. The last property shows that the expectation is of the order of $\frac{1}{J}$ whereas the variance is of the order of  $\frac{1}{J^2}$. The choice we make to preserve a stochastic part in our limit equation is to consider  the infinitesimal time variation is of the order of $\frac{1}{J^2}$ .Other choices would have led to a Piecewise Deterministic Markov process in which only the parameters $s$, $m$, $p$ brings randomness. It will be left for further study.\\

\textbf{A scale in $\frac{1}{J^2}$ and the weak selection and immigration}.\\
Let $\Delta_{t}=\frac{1}{J^2}$ et $n=tJ^2$.\\
In this case, $\E[X_{t+\Delta_{t}}-x|X_{t}=x]=O(J)$ and the limit would be infinite. To hope for a finite term and thus to observe the influence of $ s $ and $ m $ in our limit DSE, we must assume that $ s $ and $m$ are inversely proportional to $ J $.
We now assume that  migration and speciation are weak.
\begin{prop}\label{exp-var}
Let $\Delta_{t}=\frac{1}{J^2}$ and $n=tJ^2$, note
$m'(t)=m(t)\times J$ and $s'(t)=s(t)\times J $ .
So when $J$  goes to infinity:\\
\begin{enumerate}

\item $\E[X^{i}_{t+\Delta_{t}}-x^{i}|X_{t}=x]=\left(m_{t}'(p_{t}^{i}-x^{i})+x^{i}(s'^{i}_{t}-\sum\limits_{k\in \mathbb{S}}x^{k}s'^{k}_{t})\right)\times \Delta_{t}+o( \Delta_{t})$

\item $\mbox{\rm Var}[X^{i}_{t+\Delta_{t}}-x^{i}|X_{t}=x]=2x^{i}(1-x^{i})\Delta_{t} +o( \Delta_{t})$

\item $\mbox{\rm Cov}[X^{i}_{t+\Delta_{t}},X^{l}_{t+\Delta_{t}}|X_{t}=x]=-2x^{i}x^{l}\Delta_{t}+o( \Delta_{t}) \quad  \forall i \neq l.$
\end{enumerate}

\end{prop}
We will now introduce notations and assumptions. Remark that a very recent work by Bansaye et al \cite{scallim} considered a very general framework for population model convergence in random environment, but in their work the environment is usually i.i.d. whereas we are in a Markovian setting. To get shorter statement and proofs, we will make considerable simplifications for our main assumption. 

{\bf Assumption (A)}.\begin{itemize}
\item the process $(p_n)$ is assumed to be constant, which corresponds to a non evolving pool of immigration;
\item the process $(m^J)$ is an autonomous Markov chain, and consider its rescaled piecewise linear extension $\tilde m_t^J=J\,m^{J}_{\lfloor tJ^2\rfloor}$, which is assumed to take values in a finite set $E_m$ and uniformly bounded (in $J$). Let denote $P_{m^J}(m,m',t)$ its transition probabilities and assume for all $m\not=m'$,$$\lim\limits_{J\rightarrow\infty}P_{m^J}\left(m,m',\frac{1}{J^{2}}\right)\times J^{2} =Q(m,m')$$ is well defined;
\item the process $(s^J)$ is an autonomous Markov chain, and consider its rescaled piecewise linear extension $\tilde s_t^J=J\,s^{J}_{\lfloor tJ^2\rfloor}$, which is assumed to take values in a finite set $E_s$ and uniformly bounded (in $J$). Let denote  $P_{s^J}(s,s',t)$ its transition probabilities and assume for all $s\not=s'$,$$\lim\limits_{J\rightarrow\infty}P_{s^{J}}\left(s,s',\frac{1}{J^{2}}\right)\times J^{2} =Q^{s}(s,s')$$ is well defined.
\item The processes $(m_n)$ and $(s_n)$ are supposed to be independent.
\end{itemize}
These assumptions about the limits of transitions probabilities are intended to ensure the convergence  in law of $s$ and $m$ towards Markovian jump processes when $J$ goes to infinity.
  

Denote again $U^{J}_{t}=\begin{pmatrix}
X^{J}_{t} \\ s^{J}_{t}\\m_{t}^{J}
\end{pmatrix}$ taking values in  $E=E_{x}\times E_{s}\times E_{m}$  a compact set of $\mathbb{R}^{2S+1}$  and for  $\Gamma \in E$,  we define for $k\in \mathbb{N}$, $t\in [\frac{k}{J²},\frac{k+1}{J²}[$ , $\pi_{J}(U_{t},\Gamma)=\pi_{J}(U_{\frac{k}{J^{2}}},\Gamma)=\mathbb{P}(U^{J}_{\frac{k+1}{J^{2}}}\in \Gamma| U^{J}_{\frac{k}{J^2}})$.

We are now in position to state our diffusion approximation result.

\begin{theorem}
\label{th:diffapprox}
Assume {\bf (A)} 
then when $J$ goes to infinity the sequence of processes $(X_{t}^{J})$ converges in law  to the process $(X_{t})$ whose coordinates are solutions of the following stochastic differential equation 
\begin{equation}
\begin{pmatrix}
dX^{1}_{t}\\
\vdots\\
dX^{S}_{t}
\end{pmatrix}
=
\begin{pmatrix}
m_{t}\times (p_{t}^{1}-X^{1}_{t})+X^{1}_{t}\left(s_{t}^{1}-\sum\limits_{k\in \mathbb{S}}X^{k}s^{k}_{t}\right)\\
\vdots\\
m_{t}\times (p^{S}_{t}-X^{S}_{t})+X^{S}_{t}\left(s^{S}_{t}-\sum\limits_{k\in \mathbb{S}}X^{k}s^{k}_{t}\right)
\end{pmatrix}dt
+ \sigma(X_{t}) dB_{t} \label{eqprincipale2} 
\end{equation}

where $\sigma$ is such that $\sigma.\sigma^{*}=a$ with $a_{i,j}=-2x^{i}x^{j}$ if $i\neq j$ and $a_{i,i}=2x^{i}(1-x^{i})$ and where
$s_{t}=\lim s_{t}^{J}J$ and $m_{t}=\lim m_{t}^{J}J$ are  the Markovian jump processes with generators $Q^{s}$ and $Q^{m}$ and for initial conditions $s_{0}=\lim Js_{0}^{J}$ and $ m_{0}=\lim J m_{0}^{J}$. \\

\end{theorem}

\begin{proof}
The proof is in Section \ref{diffapprox} and relies on the usual Martingale Problem Method.
\end{proof}
Let us give some remarks
\begin{enumerate}
\item We may also consider the proportions $p_{i}$ random but their law would be $J$ dependent (through the change of time) and it is in disharmony with the biological model which assumes the pool independent of the community size.
\item Other types of processes for $s$ would have led to similar results, for example we could take a diffusion, with obvious modifications.
\item It is possible to give an upper bound of the error made by performing the diffusion approximation, by a direct approach not relying on the martingale problem method. It will be the purpose of \cite{ggp2018}.
\end{enumerate}

\subsubsection{Simpson index}
\quad\\
Our main object to quantify the biodiversity will be the Simpson index :
\begin{equation}
\label{simpson}
\mathscr{S}_t=\sum\limits_{i=1}^{S+1}(x_{t}^{i})^{2}
\end{equation}
with  $x_{t}^{S+1}=1-\sum\limits_{i=1}^{S}x_{t}^{i}$. \\
Notice that this expression is the limits of the discrete Simpson index when $J$ goes to the infinity. Its dynamic is given by a non autonomous stochastic differential equation.
\begin{prop}
Denote as usual  $p_{t}^{S+1}=1-\sum\limits_{i=1}^{S}p_{t}^{i},\quad \textnormal{and} \quad X^{S+1}_{t}=1-\sum\limits_{i=1}^{S}X^{i}_{t}$.\\
So the  Simpson index is solution of the following equation :\\
\vspace{2mm}
$$d\mathscr{S}_{t}=2(1-\mathscr{S}_{t})-2\sum\limits_{i=1}^{S}s_{t}^{i}X_{t}^{i}(\mathscr{S}_{t}-X_{t}^{i})+2m_{t}\left(\sum\limits_{i=1}^{S+1}p^{i}X^{i}_{t}-\mathscr{S}_{t}\right)dt+dM_{t} \label{eqSn}$$\\
where $M_{t}$ is a martingale.\\
The drift is composed of three terms, the first is the drift in the neutral case without immigration (autonomous equation), the second is a term due to the presence of selection only and the third to the presence of  immigration.
\end{prop}

This proposition follows from It\^o's calculus and details may be found in Section 5. As we will consider mainly the evolution of the expectation of the Simpson index, we do not describe the martingale term.

Now that the large limit SDE is established, we may consider the approximation of the Simpson index. We will first consider a simplified case, but which contains all the main difficulties: the two species case.

\section{Approximation of the Simpson index in the quenched or deterministic case: the two species case.}
In this part, we study  a population of only two species. The equations obtained are in dimension one and the quantities are easier to calculate. It is a basic example to understand the dynamic in greater dimension.\\
In all this section we will suppose that the {\it selection and immigration processes are deterministic}, which also amounts to consider the quenched case, i.e. we fix the random environment (immigration and selection), as our goal will be to give a numerical method  to approximate $\E[\mathscr{S}_{t}]$ at a lower cost. We will see in the next section how to consider the annealed case for a very particular case of selection and immigration.
In a second part, we will look for constant selection and immigration the behaviour of the process in long time.\\

The one dimensional Simpson index is $$\mathscr{S}_{t}=X_{t}^{2}+(1-X_{t})^{2}$$
Following the result of the previous section we are thus interested in the case where $X_{t}$ and $\mathscr{S}_{t}$ dynamics are given by

\begin{equation}
dX_{t}=m_{t}(p_{t}-X_{t})dt+s_{t}X_{t}(1-X_{t})dt+\sqrt{2X_{t}(1-X_{t}}dB_{t} \label{eq2}
\end{equation}
\begin{align}
d\mathscr{S}_{t}=&4X_{t}(1-X_{t})\times\left(1+s_{t}(X_{t}-\frac{1}{2})\right)+2m_{t}(p_{t}-X_{t})(2X_{t}-1)dt\nonumber\\  
&+2(2X_{t}-1)\, \sqrt{2X_{t}(1-X_{t})}dB_{t} \label{eqS2} 
\end{align}
where $m_t$ and $s_t$ are the rescaled limit immigration and selection processes.
\begin{remarque}
There is a first interesting feature when analysing the instantaneous behaviour of the dynamic of $\mathscr{S}_t$ in the case where there is no immigration (and thus the Simpson index will tend to 1): when $|s_t|<2$, the drift of $\mathscr{S}_t$ is always positive so that a even variability if the selection is small does not change the trend to non diversity. At the opposite, if the selection is sufficiently strong it may change locally the behaviour of the Simpson index, and we may thus imagine that change of fitness may lead to oscillation of the Simpson index. We will illustrate this phenomenon numerically.
\end{remarque}

We will now concentrate on a method to approximate all the moments of $X_t$, and thus an approximation of $\E[\mathscr{S}_t]$.

%

\subsection{Approximation of the moments of $X_{t}$}
 
\subsubsection{The approximation theorem}
\quad\\
As remarked earlier the momentum of $ \mathscr{S} _ {t} $ can not be calculated directly, as the equation of $\E[ \mathscr {S} _ {t} ]$ is not autonomous\eqref{eq2}. However we only need to evaluate the first two moments of $X_{t}$. We will see that it will be more difficult than planned. Indeed, taking expectation in  \eqref{eqS2} (recalling that $m_t$ and $s_t$ are considered deterministic), we get:

\begin{equation}\label{eqmmt}
d\E[X_{t}]=m_{t}p_{t}-m_{t}\E[X_{t}]
+s_{t}\big(\E[X_{t}]-\E[X_{t}^{2}]\big)dt
\end{equation}
By analyzing \eqref{eqmmt}, we cannot express the first momentum of $X_{t}$ without  the second moment and when considering the second moment, the third will appear and so on. It is thus impossible to express the momentum of $X_{t}$ as the solution of an autonomous equation (except for the trivial case where $s_t=0$)\\
Nevertheless, the following theorem gives a way of obtaining an approximation of the first moments of $X_{t}$ by solving a differential system whose size is all the greater that one wishes to be precise. 

\begin{theorem}
Denote
$\tilde{A}^n_t$ the tridiagonal matrix whose coefficients are given by $\tilde{a}_{i+1,i}=(i+p_tm_t)(i+1)$, $\tilde{a}_{i,i}=i(s_t-i+1-m_t)$, $\tilde{a}_{i,i+1}=-is_t$ for i in $\{1,...,n-1\}$  and $\tilde{a}_{n,n}=-n(n-1+m_t)$



Let consider the following system of ordinary differential equations
$$d\tilde M_{t}=\tilde{A}^n_{t}\times \tilde M_{t}dt+\mathscr{C}_tdt$$
 where $\mathscr{C}_t=(m_tp_t,0,...,0)^t$.\\

So for $j$ fixed, the $j^{th} $ coordinate of the solution $\tilde M_{t}$ converges when $ n $ (the size of the differential system) tends to infinity towards $j^{th}$ momentum of $X_{t}$. The error committed is at most  $$\frac{\sqrt[]{n}\|s\|_{\infty}^{n-1}}{(n-1)!}.$$
\end{theorem}

As seen from the upper bound, the convergence is quite fast and even enable to approximate the Laplace transform of $X_t$ quite efficiently. It is mainly due to the fact that we have only two species here. We will see in the next section what happens for three species and will explain how it deteriorates with the number of species.

\subsubsection{Proof of the theorem}

We will begin by considering the error between the solution of our approaching system and the real solution. For practical reasons, each coordinate of the error is multiplied by a coefficient independent of the system size. Let us begin by giving the (non autonomous) system of ordinary differential equations verified by the moments of $X_t$

\begin{prop} Let $M_{t}$ being the vector having for coordinate $i$ the $\E[X_t^i]$ (up to $n^{th}$ moment). $M_{t}$ is  solution of
$$dM_{t}=\tilde{A}_{t}^{n}\times M_{t}dt+\mathscr{C}_tdt+B_{t}dt$$
where $B_{t}=(0,...,0, ns_{t}\E[X_{t}^{n}(1-X_{t})])^T$
\end{prop}

 \begin{proof}
 \quad\\
It is of course a simple consequence of It\^o's Formula $X^{i}_{t}$:  
\begin{align*}
 dX_{t}^{i}=&iX_{t}^{i-1}\left(\left(m_{t}(p_{t}-X_{t})+s_{t}X_{t}(1-X_{t})\right)\times (i-1)(1-X_{t})\right)dt+d\mathscr{M}_{t}\\=&-is_{t}X_{t}^{i+1}+i\left(s_{t}-(i-1)-m_{t}\right)X_{t}^{i}+i(i-1+m_{t}p_{t})X_{t}^{i-1} dt+d\mathscr{M}t_{t}
 \end{align*}
 where $\mathscr{M}$ is a martingale and by taking the expectation

 $$d\E[X_{t}^{i}]=-is_{t}\E[X_{t}^{i+1}]+i(s_{t}-(i-1)-m_{t})\E[X_{t}^{i}]+i(i-1+m_{t}p_{t})\E[X_{t}^{i-1}] dt$$
  to recover the coefficients of the matrix $\tilde{A}_{t}^{n}$.
 
 \end{proof}

 We define the vector of errors  $\Delta_{t}^{n}$  by $$\Delta_{t}^{n}(i)=(s_{t})^{i-1}\times\frac{M_{t}^{i}-\tilde{M_{t}}^{i}}{(i-1)!},$$
 and we introduce also $\tilde\mu^n_t=(0,...,0,\mu_t^n)^T$ with $$\mu_{t}^{n}=\frac{(s_{t})^{n-1}}{(n-1)!}\times ns_{t}E[X_{t}^{n}(1-X_{t})].$$

Note that when we multiply, for  $i$ fixed,  the difference  $M_{t}^{i}-\tilde{M_{t}^{i}}$ by a coefficient independent of $n$, the speed of convergence of the $i^{th}$ coordinate of the error change but not the fact that this quantity tends towards $0$. On the other hand, it forces the dependent coordinates of $n$ to tend to $0$. So we just need to prove that $\|\Delta_{t}^{n}\|$ goes to  $0$, which is not the case for $M_{t}^{i}-\tilde{M_{t}^{i}}$.

 \begin{prop}
 \quad
  $\Delta_{t}^{n}$ is the solution of the following differential system:
  
  \begin{equation}
   \left\{
    \begin{aligned}
         &d\Delta_{t}^{n}=A_{t}^{n}\times \Delta_{t}^{n}dt+\tilde{\mu}_{t}^{n}dt\\
         &\Delta^{n}_{0}=0 
    \end{aligned}
  \right.
   \label{Em}
\end{equation}
  
and $$A_{t}^{n}=\begin{pmatrix}
a_{t}^{1}&b_{t}^{1}&&&(0)\\
c_{t}^{1}&a_{t}^{2}&b_{t}^{2}\\
&\ddots&\ddots&\ddots&\\
&&c_{t}^{n-2}&a_{t}^{n-1}&b_{t}^{n-1}\\
(0)&&&c_{t}^{n-1}&a_{t}^{n}
\end{pmatrix}$$

with
 \begin{align*}
 a_{t}^{i}&=i\big(s_{t}-(i-1)-m_{t}\big), \quad i< n\\
 a_{t}^{n}&=-n\big((n-1)+m_{t}\big)\\
 b_{t}^{i}&=-i^{2} \quad i< n\\
 c_{t}^{i}&=s_{t}i(1-\frac{m_{t}p_{t})}{i-1})\quad i< n\\
 \end{align*}





 \end{prop}

\begin{proof}

\quad\\

Indeed, by substraction,  $$d(M_{t}^{n}-\tilde{M}_{t}^{n})=\tilde{A_{t}^{n}}\times (M_{t}^{n}-\tilde{M}_{t}^{n})+\begin{pmatrix}
 0\\ \vdots\\ 0\\s_{t}nE[(X_{t})^{n}(1-X_{t})]\end{pmatrix}.$$
 
 This involves for all $ k<n$\\ 
 \begin{align*}
d(M_{t}^{k}-\tilde{M}_{t}^{k})&=k(k-1+p_{t}m_{t})(M_{t}^{k-1}-\tilde{M}_{t}^{k-1})\\&+k(s_{t}-(k-1)-m_{t})(M_{t}^{k}-\tilde{M}_{t}^{k})-ks_{t}(M_{t}^{k+1}-\tilde{M}_{t}^{k+1}).\\
  \end{align*}
 Next, we multiply by$\frac{(s_{t})^{k-1}}{(k-1)!}$,
 
 \begin{align*}
 d(\Delta_{t}^{k,n})&=\frac{(M_{t}^{k-1}-\tilde{M}_{t}^{k-1})s_{t}^{k-2}}{(k-2)!}\frac{s_{t}k(k-1+p_{t}m_{t})}{k-1}\\&\qquad+\frac{ks_{t}^{k-1}}{(k-1)!}\left(s_{t}-(k-1)-m_{t})(M_{t}^{k}-\tilde{M}_{t}^{k}\right)\\&\qquad
 -\frac{s_{t}^{k}(M_{t}^{k+1}-\tilde{M}_{t}^{k+1})}{(k)!}k^{2}\\
 &=ks_{t}\big(1-\frac{m_{t}p_{t})}{k-1}\big)\Delta_{t}^{k-1,n}+k(s_{t}-(k-1)-m_{t})\Delta_{t}^{k,n}-k^{2}\Delta_{t}^{k+1,n}.
 \end{align*}

 Now if k=n,
 \begin{eqnarray*}
 d\Delta_{t}^{n,n}&=&n\left[s_{t}\left(1-\frac{m_{t}p_{t}}{n-1}\right)\Delta_{t}^{n-1,n}-(m_{t}+n-1)\Delta_{t}^{n,n}\right.\\
 &&\qquad\left.+\frac{s_{t}^{i-1}}{(n-1)!} s_{t}\E[X_{t}^{n}(1-X_{t})]\right]
 \end{eqnarray*}
 
We thus find the coefficients of the previous equation system.
\end{proof}

We may now provide the solution to the system of equation for the error.

\begin{prop}

The solution of $\eqref{Em}$ can be written $$\Delta_{t}^{n}=\int_{0}^{t}\exp\left(\int_{q}^{t}A_{u}^{n}\right)\tilde{\mu}_{q}^{n}dq$$ 
\end{prop}
\begin{proof}
First, the solution of $\eqref{Em}$ can be written as the solution of the homogeneous system and a particular solution.
$$\Delta_{t}^{n}=\mathcal{K}\exp\left(\int_{q}^{t}A_{u}^{n}du\right)+\exp\left(\int_{0}^{t}A_{u}^{n}du\right)\int_{0}^{t}\exp\left(-\int_{0}^{q}A_{u}^{n}du\right)\tilde{\mu}_{t}^{n}(q)dq
$$ As $\Delta_{0}^{n}=0$ necessarily  $\mathcal{K}=0$ and $\Delta_{t}^{n}=\int_{0}^{t}\exp\left(\int_{q}^{t}A_{u}^{n}du\right)\tilde{\mu}_{t}^{n}(q)dq.$
\end{proof}


We will now show that, for a fixed time interval,  $\|\Delta_{t}^{n}\|_{2}$ is uniformly bounded by a quantity which goes to 0 when n goes to infinity. In the following, $\|.\|_{F}$ stands for the Frobenius norm. Thanks to the formula of the previous proposition:

\begin{align*}
\|\Delta_{t}^{n}(t)\|_{2}^{2}&=\left\|\int_{0}^{t}\exp\left(\int_{q}^{t}A_{u}^{n}du\right)\tilde{\mu}_{q}^{n}dq\right\|_{2}^{2}\\
&=\sum\limits_{k=1}^{n}\left(\int_{0}^{t}(\exp\left(\int_{q}^{t}A_{u}^{n}du\right)\tilde{\mu}_{q}^{n})dq\right)_k^{2}\\
&\leqslant  t\times \sum\limits_{k=1}^{n}\int_{0}^{t}\left(\exp\left(\int_{q}^{t}A_{u}^{n}du\right)\tilde{\mu}_{q}^{n}\right)_{k}^{2}dq\\
&\leqslant t  \int_{0}^{t}\left\|\exp\left(\int_{q}^{t}A_{u}^{n}du\right)\tilde{\mu}_{q}^{n}\right\|_{2}^{2}dq        \\
&\leqslant t \int_{0}^{t}\left\|\exp\left(\int_{q}^{t}A_{u}^{n}du\right)\right\|_{F}^{2}\,\|\tilde{\mu}_{q}^{n}\|_{2}^{2}
dq\\
&\leqslant t^{2} \times \sup\limits_{x\in[0,t]}\left\|\exp\left(\int_{x}^{t}A_{u}^{n}du\right)\right\|_{F}^{2}\times \sup\limits_{x\in[0,t]}|\mu_{x}^{n}|^{2}
\end{align*}
 
The following easy lemma gives us an upper bound for the two previous norms:
 
 \begin{lemma}
\quad\\
If $X \in [0,1]$ then $sX^{n}(1-X)=O(\frac{1}{n})$ and there is a constant $ c_ {1} $ independent of n such that$\sup\limits_{x\in[0,t]}|\mu_{x}^{n}|\leqslant \sup\limits_x\in[0,t]s_{x}^{n-1}\frac{c_{1}}
{(n-1)!}$
\end{lemma}
   
\begin{proof}
The maximum of $X^{n}(1-X)$ on $[0,1]$ is achieved in $\frac{n}{n+1}$ and it's worth   $(1-\frac{1}{n+1})^{n}\times \frac{1}{n+1}$.\\
This quantity is of the order of $\frac{1}{n}$ when $n$ goes to infinity.\\
The upper bound of  $\mu_{x}^{n}$ follow.
\end{proof}

Next,  \begin{align*}
\left\|\exp\left(\int_{x}^{t}A_{u}^{n}du\right)\right\|_{F}^{2}=&Tr\left(\exp\left(\int_{x}^{t}A_{u}^{n}du\right)^{T}\times \exp\left(\int_{x}^{t}A_{u}^{n}du\right)\right)\\
 =&Tr\left(\exp\left(\int_{x}^{t}A_{u}^{n}+(A_{u}^{n})^{T}du\right)\right)\\
 =&\sum\limits_{\lambda_{i} } \exp\left(\lambda_{i}^{x}\right)\\
 \leqslant & n\times \max\limits_{i \in \{1...n\}}(\
 \exp(\lambda_{i}^{x}))
\end{align*}
where the $\lambda_{i}^{x}$ are the eigenvalues of  $\int_{x}^{t}A_{u}^{n}+(A_{u}^{n})^{T}du$ .\\

Then if $\max\limits_{i \in \{1...n\}}\exp(\lambda_{i}^{x})$ is independent of $n$, there is a constant $c_{2}$ independent of $n$ verifying $$ \sup\limits_{x\in[0,t]}\left\|\exp(\int_{x}^{t}A_{u}^{n}du)\right\|_{F}\leqslant \sqrt[]{n}\times c_{2}$$

And then there exist a constant $C$ such as  $$\|\Delta_{t}^{n}(t)\|_{2}\leqslant\sqrt[]{n} c_{2}\times \frac{c_{1}\|s\|_{\infty}^{n-1}}{(n-1)!}\leqslant\frac{C\sqrt[]{n}\|s\|_{\infty}^{n-1}}{(n-1)!} $$

It remains to show that the eigenvalues of $\int_{x}^{t}A^{n}_{u}+(A_{u}^{n})^{T}du$ (which are real) have an upper bound  independent of n.

\begin{prop}

Let $A_{t}^{n}+(A_{t}^{n})^{T}=
\begin{pmatrix}
a_{t}^{1}&b_{t}^{1}&&&\\
b_{t}^{1}&a_{t}^{2}&b_{t}^{2}\\
&\ddots&\ddots&\ddots&\\
&&b_{t}^{n-2}&a_{t}^{n-1}&b_{t}^{n-1}\\
&&&b_{t}^{n-1}&a_{t}^{n}
\end{pmatrix}
$
 with
 \begin{align*}
 a_{t}^{i}&=2i\big(s_{t}-(i-1)-m_{t}\big), \quad i< n\\
 a_{t}^{n}&=-n\big((n-1)+m_{t}\big)\\
 b_{t}^{i}&=s_{t}i(1-\frac{m_{t}p_{t})}{i-1})-(i-1)^{2}
 \end{align*}

 then the eigenvalues of $\int_{x}^{t}(A_{u}^{n})^{T}+A_{u}^{n}du$ are uniformly bounded in $n$ .
 \end{prop}
 
 \begin{proof} 
 To prove this result we need to use Gershgorin's disk. The eigenvalues of  $\int_{x}^{t}(A_{u}^{n})^{T}+A_{u}^{n}du$ are included in the union of the disk $D_{i}$ whose centers are the  $i^{th}$ term on the diagonal ($\int_{x}^{t}a_{u}^{i}du$) and for radius the sums of the coefficients  norms on the $i^{th}$ line except the diagonal term$(|\int_{x}^{t}b_{u}^{i-1}du|+|\int_{x}^{t}b_{u}^{i}du|)$. It is important to consider the forms of the discs in our case.
As we have to show that the eigenvalues have a upper bound independently of n. We just need to look at the shape of the discs for $i $ and $n$ big enough. From the matrix if $i<n$, their centers are  
$$2i\left(\int_{x}^{t}s_{u}du-(i-1)(t-x)-\int_{x}^{t}m_{u}du\right)$$ and their radius are equal to 
\begin{align*}
\left|\int_{x}^{t}b_{u}^{i-1}du\right|+\left|\int_{x}^{t}b_{u}^{i}\right|&=\left|\int_{x}^{t}s_{u}i\left(1-\frac{m_{u}p_{u}}{i-1}\right)-(i-1)^{2}du\right|\\&\qquad+\left|\int_{x}^{t}s_{u}(i+1)\left(1-\frac{m_{u}p_{u}}{i}\right)-i^{2}du\right|\\
&=\left|i\int_{x}^{t}s_{u}du-(i^2-2i)(t-x)+i\int_{x}^{t}m_{u}du\right|\\&\qquad+\left|i\int_{x}^{t}s_{u}du-i^2(t-x)-i\int_{x}^{t}m_{u}du\right|+O(1)\\
&=2i\left((i-1)(t-x)-\int_{x}^{t}s_{u}du+\int_{x}^{t}m_{u}du\right)+O(1).
\end{align*}
So the maximum value for an eigenvalue of $\int_{x}^{t}(A_{u}^{n})^{T}+A_{u}^{n}du$ belonging to $D_{i}$ is
\begin{align*}
&\int_{x}^{t}a_{u}^{i}du+\left|\int_{x}^{t}b_{u}^{i-1}du\right|+\left|\int_{x}^{t}b_{u}^{i}du\right|\\&\qquad\qquad=2i\left(\int_{x}^{t}s_{u}du-(i-1)(t-x)-\int_{x}^{t}m_{u}du\right)\\&\qquad\qquad\qquad+2i\left((i-1)(t-x)-\int_{x}^{t}s_{u}du+\int_{x}^{t}m_{u}du\right)+O(1)\\ &\qquad\qquad=O(1)
\end{align*} as soon as $i$ is big enough. If $i=n$ the same reasoning is still working. So the eigenvalues of  $(A^{n})^{T}+A^{n}$ have an upper bound independent of  $n$. 
  \end{proof}
 
This property concludes the proof and we get the following upper bound$$\|\Delta_{t}^{n}(t)\|_{\infty}\leqslant\|\Delta_{t}^{n}(t)\|_{2}\leqslant\frac{C\sqrt[]{n}\|s\|_{\infty}^{n-1}}{(n-1)!} $$

The constant $ C $ depends on time (exponentially) and therefore this algorithm will be less accurate if we look at the behavior of the process in long time. This result gives a satisfactory approximation of the $ \mathscr {S} _ {t} $'s moment. Convergence is very fast, and the algorithm boils down to solving a linear differential system. To ensure the interest of this method we can compare the expectation of Simpson's index obtained by a Monte Carlo method from the discrete model to that obtained with this approximation. The Figure 1 presents such an approximation.

\begin{figure}[htbp]
\begin{minipage}[c]{.45\linewidth}
\begin{center}
\includegraphics[scale=0.35]{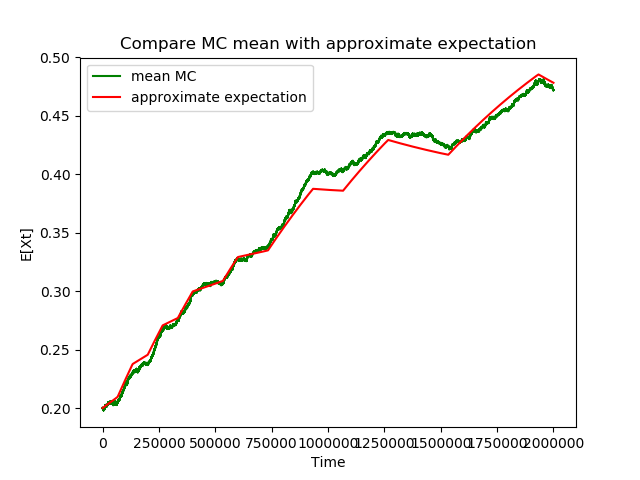}
\end{center}
\end{minipage}
\hfill
\begin{minipage}[c]{.45\linewidth}
\begin{center}
\includegraphics[scale=0.35]{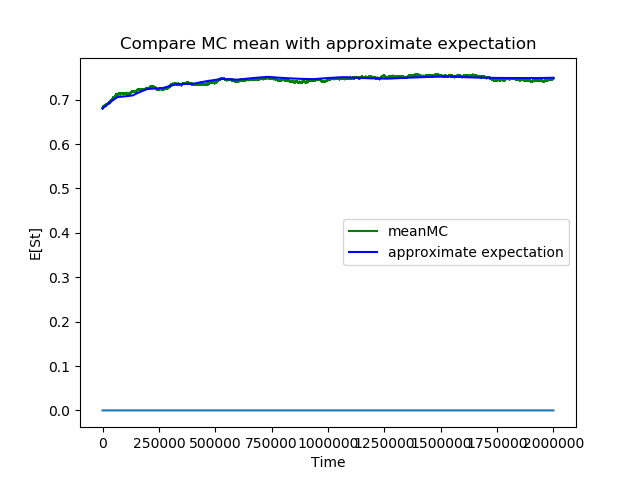}
\end{center}
\end{minipage}
\label{fig:image1}
\caption{Are plotted the approximate values of $E[\mathscr{S}_{t}]$ and $E[X_{t}]$ by the precedent method from the approximation in large population and by MC method from the discreet model. The number of simulated trajectories for  MC mean is $ 500$ (red and blue), $J=1000$,$m=2$,$p=0.5$,$X_{0}=0.2$, $s$ switches between $2$ and $-2 $ at regular time intervals, the size of the approaching linear system is $100$.}
\end{figure}

\subsection{Numerical applications}
\quad\\
The simulations presented in this part are obtained from the previous theorem.
The values of $ s, m, p $  are those of the large population approximation and not that of the discrete model.
The size of the approaching system will be usually between 80 and 144 depending on the needed precision. 
\subsubsection{Influence of $s$ on Simpson Index.}
\quad\\
In this part, $p=0.5$. Now  we know how to approximate the expectation of $ \mathscr {S} _ {t} $, so we can check the influence of  $ s $ on this quantity. Let us make precise a statement enounced when deriving $\mathscr{S}_t$.
\begin{prop}
\quad\\
If $\|s\|_{\infty}$ is smaller than $2$ and if $m=0$,  $E[ \mathscr {S} _ {t} ]$ is increasing.
\end{prop}
\begin{proof} 
\quad\\
If we refer to the equation of $ dE [\mathscr {S} _ {t}] $ (cf \eqref {eqS2}), we see that if $\|s\|_{\infty}<2$, the quantity $ 4X_{t}(1-X_{t})\times\left(1+s_{t}(X_{t}-\frac{1}{2})\right)$ is positive  whatever the initial condition, so Simpson's index mean is always  growing .
\end{proof}
In other words, selection alone can not bring about a renewal of biodiversity.

On the other hand if for some $t$, $s_{t}> 2$, $ E [\mathscr {S} _ {t}] $ can decrease.
 In this case the more the selection is important, more the decay is pronounced.
 We will see in the last part that this phenomenon can be generalized  to a larger number of species. Figure 2 present these different behaviours with respect to $s$ and Figure 3 the combination of initial parameters and selection under which the Simpson's index is decreasing.

\begin{figure}[htbp]
\begin{minipage}[c]{.45\linewidth}
\begin{center}
\includegraphics[scale=0.35]{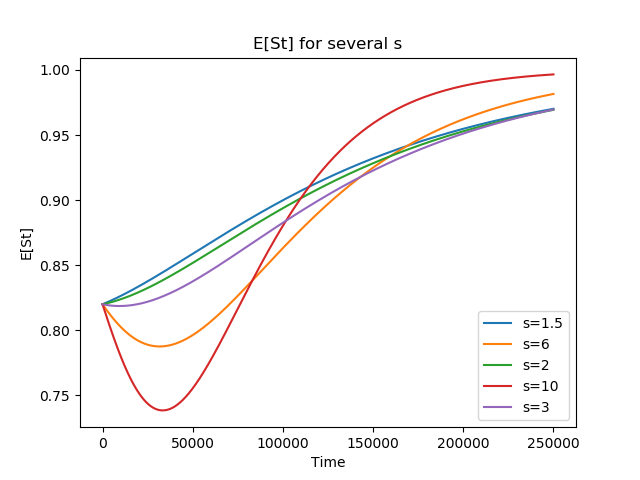}
\label{fig:image2}
\caption{Several trajectories of $E [\mathscr {S} _ {t}]$ are drawn for different s. $X_{0}=0.1$, m=0.  }
\end{center}
\end{minipage}
\hfill
\begin{minipage}[c]{.45\linewidth}
\begin{center}
\includegraphics[scale=0.35]{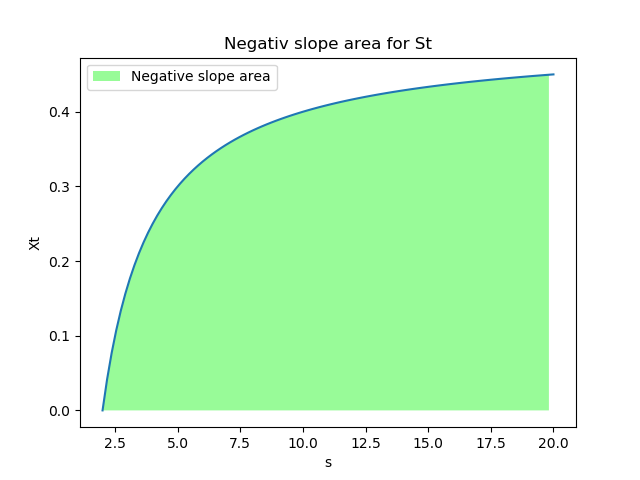}
\label{fig:image3}
\caption{The colored area represents the pairs $(X, s)$ for which the slope of the Simpson's index will be negative.}
\end{center}
\end{minipage}
\end{figure}

\subsubsection{Approximation of $T_{1}$, $T_{0}$, $T_{1,0}$.}
\quad\\
In the special case where $m_{t}=0$, a species inevitably invaded the community in a finite time. We define by $T_{1}$, (respectively $T_{0}$),  the smallest time from which the process $X_{t}$ reaches $1$ (respectively $0$) and $T_{1,0}=min( T_{1}, T_{0}) $.\\
Thanks to the approximation of moments we can  obtain an approximation of the  $ T_ {1} $ distribution function .
In fact we know $X_{t}^{n}$ for $n$ big enough and so since $\lim\limits_{n\rightarrow \infty}E[(X_{t})^{n}]=1\times \pp(T_{1}<t)$ we obtain $P(T_{1}<t)$.
The same way with $1-X_{t}$ gives $\pp(T_{0}<t)$. Figure 4 gives an approximation of the distribution function of $T_1$ when there is no immigration.

\begin{figure}[htbp]
\begin{minipage}[c]{.45\linewidth}
\begin{center}
\includegraphics[scale=0.35]{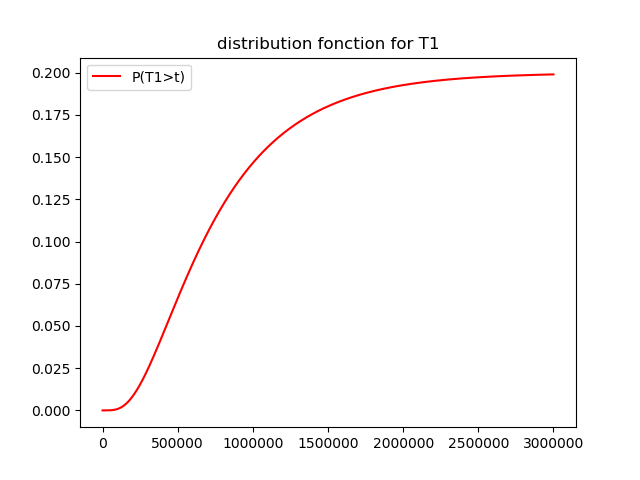}
\end{center}
\end{minipage}
\caption{Distribution function of $T_{1}$ for m=0,p=0.5,$X_{0}=0.2$, $s=2$. The size of the approaching linear system is 100}
\label{fig:image4}
\end{figure}

We can use the same method to obtain $T_{\mathscr{S}_{t}}$ the probability that $\mathscr{S}_{t}$ is equal to $1$ at time $t$.

\subsection{Long time behavior}
\quad\\
Two cases are distinguished in this part, the case $ m = 0 $ and the case $ m \neq 0 $. In this first case there is no immigration and necessarily a species invade the community. Invasion times and the probability that the species with a selective advantage will invade the community are calculated. In the second case the system admits an invariant measure, we  explain it and we specify the speeds of convergence towards this measure.
More details about the behavior of stochastic processes in long time can be found in \cite {evo} and \cite {SMB}.

\subsubsection{The case without immigration ($m=0$): absorption}
\quad\\
The results in this section are partially well known and we include them only to get a full picture of the behavior of the Simpson index.\\ Recall that $\mathscr{S}_{t}$ satisfies the equation 
\begin{equation}
d\mathscr{S}_{t}=4X_{t}(1-X_{t})\times\big(1+s_{t}(X_{t}-\frac{1}{2})\big)dt+dM_{t}.
\label{eqS2b}
\end{equation}

If there is no more immigration the states  $0$ and $1$ are absorbing, and it is then well known that the process reaches them in finite time almost surely. For more details about this refer to \cite{evo}. Let $T_{1}$ et $T_{0}$ the hit times of 1 and 0 for the random variable$X_{t}$ et $T_{1}\wedge T_{0}=T_{1, 0}$.
\begin{prop}\label{prop:absorb}
Suppose $s_{t}=s \in \mathbb{R}$,
\begin{enumerate}
\item $T_{1, 0}<\infty$ almost surely, for all initial condition $X_{0}$.\\

\item Let g be the solution of
\begin{equation}
x(1-x)g"(x)+sx(1-x)g'(x)=-1 \quad et \quad  
g(0)=g(1)=0 \label{g}
\end{equation}
then $ E_{X_{0}}[T_{1, 0}]=g(X_{0})$. \\

\item $P(T_{1}<T_{0})=\frac{e^{-sX_{0}}-1}{e^{-s}-1}.$
\end{enumerate}

\end{prop}

The proof is given in section \ref{absorb}


Let us consider some particular case which illustrates that the same behavior may be obtained with varying selection. Suppose $T>0$ and $s_{t}$ is a constant function on intervals $[kT,(k+1)T]$, $k\in \mathbb{N}$ which can take the values $ s_ {0} $ or $ -s_ {0} $ for $ s_ {0}> 0 $, randomly. We will establish a result similar to the constant case. Let us begin by the following lemma which asserts that without selection one may reach the boundary at any time.

\begin{lemma}
Consider the following process 
\begin{equation}
dX_{t}=\sqrt[]{2X_{t}(1-X_{t})}dB_{t} \label{eq2neutre}
\end{equation}
Note $T_{1}=\inf\{t,X_{t}=1\} $, an initial condition $x \in ]0,1]$ and a time $t$ .\\
Then $\mathbb{P}_{x}(T_{1}<t)>0$, in other words, $ 1 $ is accessible for $ X_ {t} $ from any non-zero initial condition and in a as little time as one wants.
\end{lemma}

\begin{proof}
Assume that  $t_{0}=\inf\{t,\mathbb{P}_{x}(T_{1}<t)>0\}>0$.\\
 $t_{0}$ is well defined. Remark that $ \E [X_ {t}] $ is constant in time because $ X_ {t} $ is a bounded martingale. Then $0<x=\lim\limits_{t\rightarrow\infty}E[X_{t}]=1\times \mathbb{P}_{x}(T_{1}<T_{0})+0$ as $\E_{x}[T_{1,0}]<\infty$  and so  $\mathbb{P}_{x}(T_{1}<\infty)>0$. Then, let  $\Delta_{t}>0$ be such that $t_{0}-2\Delta_{t}>0$. We will show then  that there is $ y $ such that$\mathbb{P}_{y}(T_{1}<\Delta_{t})>0$.\\

If that was not the case then $\forall y \in [0,1],\mathbb{P}_{y}(T_{1}<\Delta_{t})=0 $  and 
\begin{align*}
\mathbb{P}_{x}(T_{1}<t_{0})=&\E_{x}[\mathds{1}_{T_{1}<t_{0}}]=\E_{x}[\mathds{1}_{T_{1}<t_{0}}\times(\mathds{1}_{T_{1}<\Delta_{t}}+\mathds{1}_{T_{1}>\Delta_{t}})]\\
=&\E_{x}[\mathds{1}_{T_{1}<\Delta_{t}}]+\E_{x}[\mathds{1}_{T_{1}<t_{0}}\times\mathds{1}_{T_{1}>\Delta_{t}}]\\
=&\E_{x}[\mathds{1}_{T_{1}<t_{0}}\times\mathds{1}_{T_{1}>\Delta_{t}}]=\E_{x}[\mathds{1}_{T_{1}>\Delta_{t}}E[\mathds{1}_{T_{1}<t_{0}}|X_{\Delta_{t}}]]\\
=&\E_{x}[\mathds{1}_{T_{1}>\Delta_{t}}\times \mathbb{P}_{X_{\Delta_{t}}}(T_{1}<t_{0}-\Delta_{t})]\\
=&\E_{x}[\mathds{1}_{T_{1}>\Delta_{t}}  \E_{X_{\Delta_{t}}}[\mathds{1}_{T_{1}>\Delta_{t}} \cdots \E_{X_{n\Delta_{t}}}[\mathds{1}_{T_{1}>\Delta_{t}}\\&\qquad\qquad\qquad\times \mathbb{P}_{(n+1)\Delta_{t}}(T_{1}<t_{0}-(n+1)\Delta_{t})]]].
\end{align*}

Let us choose $n$ such that $t_{0}-(n+1)<\Delta_{t}$. Then $\mathds{1}_{T_{1}>\Delta_{t}}\times \mathbb{P}_{(n+1)\Delta_{t}}(T_{1}<t_{0}-(n+1)\Delta_{t})\leqslant \mathds{1}_{T_{1}>\Delta_{t}}\mathbb{P}_{(n+1)\Delta_{t}}(T_{1}<\Delta_{t})=0$  and $\mathbb{P}_{x}(T_{1}<t_{0}) =0$ which is contrary to the assumptions. Now,  we show that $\mathbb{P}_{x}(T_{1}<t_{0}-\Delta_{t})>0$:

\begin{align*}
\mathbb{P}_{x}(T_{1}<t_{0}-\Delta_{t})&=\E_{x}[\mathds{1}_{T_{1}<t_{0}-\Delta_{t}}\times(\mathds{1}_{T_{y}<t_{0}-2\Delta_{t}}+\mathds{1}_{T_{y}>t_{0}-2\Delta_{t}})]\\
& \geqslant \E_{x}[\mathds{1}_{T_{1}<t_{0}-\Delta_{t}}\times\mathds{1}_{T_{y}<t_{0}-2\Delta_{t}}]\\
&\geqslant \E_{x}[\mathds{1}_{T_{y}<t_{0}-2\Delta_{t}}E[\mathds{1}_{T_{1}<t_{0}-\Delta_{t}}|T_{y}]]\\
&\geqslant  \E_{x}[\mathds{1}_{T_{y}<t_{0}-2\Delta_{t}}\mathbb{P}_{y}(T_{1}<t_{0}-\Delta_{t}-T_{y})]\\
&\geqslant \E_{x}[\mathds{1}_{T_{y}<t_{0}-2\Delta_{t}}\mathbb{P}_{y}(T_{1}<\Delta_{t})]\\
&\geqslant  \mathbb{P}_{x}(T_{y}<t_{0}-2\Delta_{t})\times \mathbb{P}_{y}(T_{1}<\Delta_{t}).
\end{align*}
But we know that $ \mathbb{P}_{y}(T_{1}<\Delta_{t})>0$ by  the previous calculation and the local uniform ellipticity of diffusion also ensures us that $\mathbb{P}_{x}(T_{y}<t_{0}-2\Delta_{t})>0$. So we obtain $ \mathbb{P}_{x}(T_{1}<t_{0}-\Delta_{t})>0$ which is contrary to the fact that $t_{0}>0$.
So $\forall t>0, \mathbb{P}_{x}(T_{1}<t)>0$ which concludes the proof of the lemma.
\end{proof}

Of course, this result may be adapted for the process $dX_{t}=sX_{t}(1-X_{t})dt+\sqrt[]{2X_{t}(1-X_{t})}dB_{t}.$ Indeed if $s>0$ the drift goes in the right direction. Else, if $s<0$, we obtain a symmetric result by replacing in the previous reasoning $ T_ {1} $ by $ T_ {0} $.

\begin{prop}
\quad\\
Let $s_{0}>0$, $T>0$ and  $s_{t}$ a constant function on the intervals $[kT,(k+1)T]$, $k\in \mathbb{N}$ which can take the values  $ s_ {0} $ or $-s_ {0} $ values randomly.\\
Let consider the process 
\begin{equation}
dX_{t}=s_{t}X_{t}(1-X_{t})dt+\sqrt[]{2X_{t}(1-X_{t})}dB_{t} \label{eq2b}
\end{equation}
Finally  $T_{1,0}=\inf\{t,X_{t}=0 \quad ou \quad X_{t}=1\}$ then $\mathbb{P}_{x}(T_{1,0}<\infty)=1$.

\end{prop}

\begin{proof}
The idea is to show that for each time interval of size $ T $, the probability of reaching $ 1 $ or $ 0 $ is non-zero and independent of the position where the process is located. So we compare the probability that our process reaches $ 1 $ or $ 0 $ to a geometric law. So let us first show that $\forall t \in ]0,T],\exists \alpha>0$, such that $ \forall x \in [0,1], \mathbb{P}_{x}(T_{1,0}<t)>\alpha$. Suppose $s>0$ on $[0,T]$ and denote $g_{i}(x)= \mathbb{P}_{x}(T_{i}<t)$, $i \in  \{0,1\}$. Both functions are continuous,  $g_{0}$ is decreasing and $g_{0}(0)=1$,  $g_{0}(1)=0$, whereas $g_{1}$ is increasing and $g_{1}(1)=1$,  $g_{1}(0)=0$.
Then there exist a $x_{0}  \in ]0,1[$ such as $g_{0}(x_{0})=g_{1}(x_{0})=\alpha_{+}$.  And by the previous lemma since $ s $ does not vary on $] 0, T] $ we have that $\alpha_{+}>0$. A symmetric reasoning for $ s <$ 0 guarantees us the existence of a $\alpha_{-}>0$.\\
Let $\alpha =2\min(\alpha_{-},\alpha_{+})$, $\alpha$ is then strictly positive and  $ \forall x \in [0,1]$, $\mathbb{P}_{x}(T_{1,0}<t)= \mathbb{P}_{x}(T_{1}<t)+ \mathbb{P}_{x}(T_{0}<t)\geqslant \alpha$. Then, for  $t\in [nT,(n+1)T]$, using previous inequality:
\begin{align*}
\mathbb{P}_{x}(T_{1,0}<t)\geqslant& \mathbb{P}_{x}(T_{1,0}<nT)\\
\geqslant &\E_{x}[\mathds{1}_{T_{1,0}<(n-1)T}]+\E_{x}[\mathds{1}_{T_{1,0}>(n-1)T}\mathbb{P}_{X_{(n-1)T}}(T_{1,0}<T)]\\
 \geqslant&\mathbb{P}_{x}(T_{1,0}<(n-1)T)+\alpha \times \E_{x}[1-\mathds{1}_{T_{0,1}<(n-1)T}]\\
  \geqslant& (1-\alpha)\times \mathbb{P}_{X_{(n-1)T}}(T_{1,0}<(n-1)T)+\alpha\\
   \geqslant& \alpha+\alpha(1-\alpha)+....+\alpha(1-\alpha)^{n-1}
   \geqslant 1-(1-\alpha)^{n}
\end{align*}
We obtain when $n$ goes to infinity $\mathbb{P}_{x}(T_{1,0}<\infty)=1$.
\end{proof}

Even with frequent changes of fitness, a species always ends up invading the community if $m$ is zero. We may then consider the process with immigration. Note that while preparing this paper, comparable (and even more general) results were obtained (in the multi-allelic case) by Coron et al \cite{cmv}.

 \subsubsection{The case with immigration ($m\neq 0$): invariant measure}  
  \quad\\
Assume now $s$, $m$ and $p$ are constants. The long time behavior for varying selection and immigration is far more complicated and may lead to interesting behavior that will be considered in another paper. We thus consider the following process:
$$dX_{t}=m(p-X_{t})+sX_{t}(1-X_{t})dt+\sqrt[]{2X_{t}(1-X_{t})}dB_{t}$$
and will consider the long time behavior in a quantitative way, i.e. not using Meyn-Tweedie's theory, but rather via a Poincar\'e inequality.\\
Our process $ X_ {t} $ is Markovian and evolves in a range bounded by 0 and 1. At first we will ask ourselves what is the behaviour of our process in the neighbourhood of 0 and 1, by considering the criterion given by Feller cf\cite{Feller1,Feller2}. According to the values of $ m $, $ p $ and $ s $, our process will have different behaviours in the neighbourhood of $1$ and $0$.
We have already seen that if $m=0$ then $0$ and $1$ are absorbing states  reached by the process in a finite time almost surely.(The same hold if $ m $ is non-zero and if $ p $ is $0$ or $1$.) Now if $ m $ and $ p $ are not trivial, $0$ and $1$ are no longer absorbing. In other words, immigration prevents the invasion of the community by a species. Moreover, for some values of $ m $ and $ p $ these two states are not accessible, i.e the process can not reach them in a finite time.

\begin{prop}\label{prop:access}
The state $ 1 $ (respectively $ 0 $) is accessible by the process $ X_ {t} $ if and only if $ m (1-p) <1 $ (respectively $ mp <1 $) and regular otherwise.
\end{prop}
The proof will be given in section \ref{access}. In the case of inaccessible or reflective boundaries (which is our the case), the law of the process $ X_ {t} $ admit a density and converge in long time to an invariant measure. This measure has a density, denoted $ \pi $. In addition $ \pi $ is a solution of the Fokker-Planck equation:
$$0=\frac{\partial }{\partial t}\pi=-\frac{\partial\pi }{\partial y}\big(m(p-y)+sy(1-y)\big)+\frac{\partial² \pi}{\partial^{2} y}y(1-y)$$
The solution of this equation is:
\begin{equation}
\pi(y)=c\times y^{mp-1}\times (1-y)^{m(1-p)-1}\times \exp(sy) \label{mesinv2}
\end{equation}
The constant $c$ is chosen so that  $\int_{0}^{1}\pi(y)dy=1$.

The following Figures 5 and 6 show the influence of the parameters on the expectation and the variance of Simpson's equilibrium index.

\begin{figure}[htbp]
\begin{minipage}[c]{.45\linewidth}
\begin{center}
\includegraphics[scale=0.35]{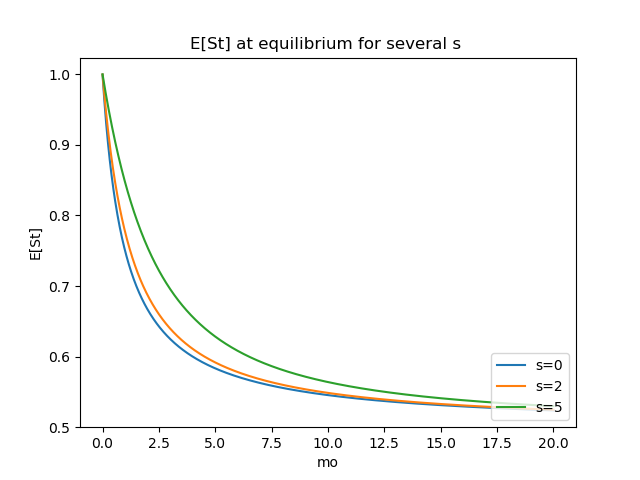}

\label{fig:image2b}
\end{center}
\end{minipage}
\hfill
\begin{minipage}[c]{.45\linewidth}
\begin{center}
\includegraphics[scale=0.35]{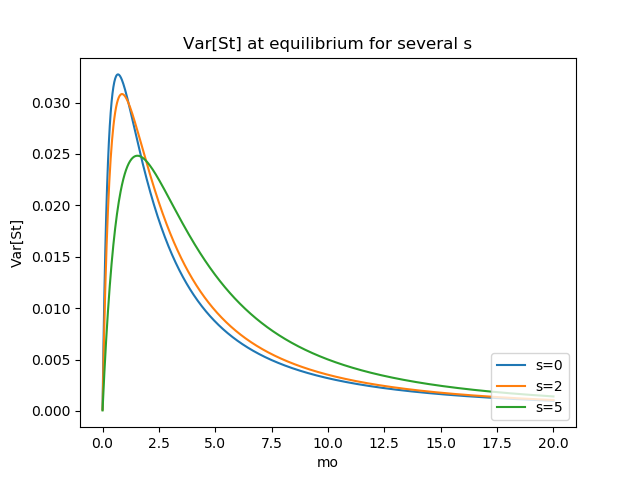}

\label{fig:image2c}
\end{center}
\end{minipage}
\caption{Expectation and variance of the Simpson index at equilibrium against m for several values of s,
p=0.5}
\end{figure}
\begin{figure}
\begin{minipage}[c]{.45\linewidth}
\begin{center}
\includegraphics[scale=0.35]{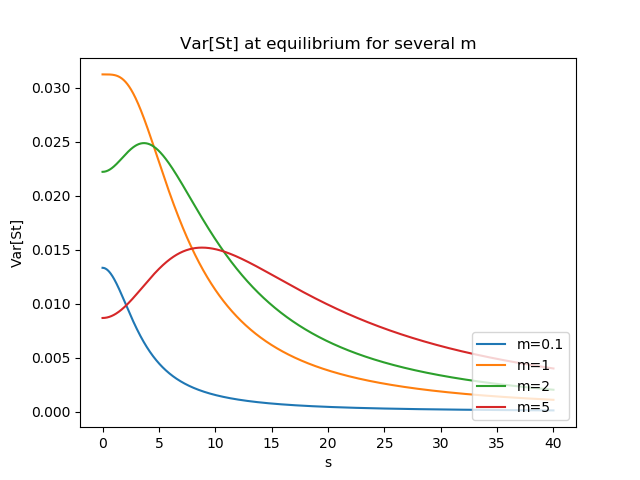}
\end{center}
\end{minipage}
\hfill
\begin{minipage}[c]{.45\linewidth}
\begin{center}
\includegraphics[scale=0.35]{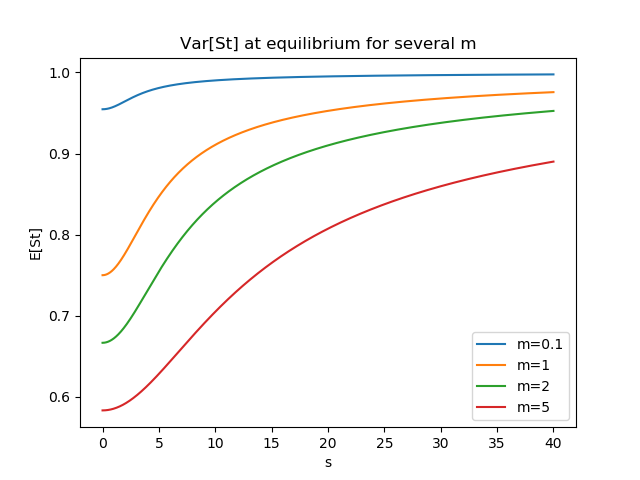}
\end{center}
\end{minipage}
\label{fig:image5}
\caption{Here are plotted the variance and the expectation of the simpson index at equilibrium against s for several values of m.
p=0.5, size of linear approaching system is 100. }
\end{figure}
\begin{figure}
\begin{minipage}[c]{.45\linewidth}
\begin{center}
\includegraphics[scale=0.35]{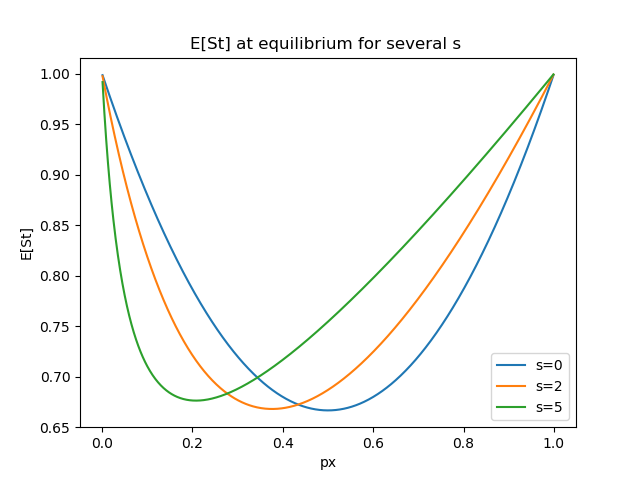}
\end{center}
\end{minipage}
\hfill
\begin{minipage}[c]{.45\linewidth}
\begin{center}
\includegraphics[scale=0.35]{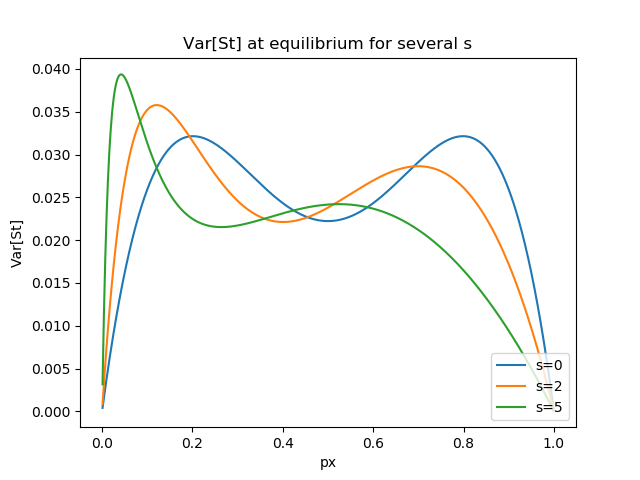}
\end{center}
\end{minipage}
\label{fig:image6}
\caption{Here are plotted the variance and the expectation of the Simpson index at equilibrium against p for several values of s.
m=2, size of linear approaching system is 100. }
\end{figure}
\vspace{30mm}
Let us now quantify the convergence to equilibrium. Recall at first that the process $ X_ {t} $ has for generator $ \mathscr {L} $ and for invariant measure $ \pi $. Denote $P_{t}f(x)=\mathbb{E}[f(X_{t})|f(0)=x]$ the associated semigroup. In fact, when $s=0$, the full spectrum is known, see for example Shimakura \cite{shimakura1977} which provides a spectral gap value $m$. It will imply an exponential convergence to equilibrium in $ L_{\pi}^{2} $. 
\begin{prop}
Let us suppose that $s>0$. The following Poincar\'e inequality is valid, i.e. for every smooth function $f$
$$\mbox{Var}_{\pi}(f)\leqslant \min\left(\frac{e^s}{m},\frac{8e^{(1-M)s}}{m}\right)\int_{0}^{1}f'^{2}(x)\times x(1-x)d\pi(x)$$
where $M$ is a median of $\pi$. As a consequence $\pi_{t}$ converge to $\pi$ in $L^{2}_{\pi}$ exponentially:
$$\mbox{Var}_{\pi}(P_{t}f)\leqslant e^{-2\max\left(\frac{e^s}{m},\frac{8e^{(1-M)s}}{m}\right)}\mbox{Var}_{\pi}(f).$$
\end{prop}

\begin{proof}
Using usual Holley-Stroock's perturbation argument we easily deduce, that the Poincar\'e constant is at most $e^s/m$. We will now use Hardy's type condition (see for example \cite{ABC}) for Poincar\'e condition, that we recall now


\begin{lemma}
Let $\pi$ and $\nu$  be two measures and $M$ the median of  $\pi$.\\
Let $$B^{+}_{M}=\sup\limits_{x>M}\int_{x}^{1}d\pi\int_{M}^{x}\frac{1}{\nu(t)}dt,\qquad B^{-}_{M}=\sup\limits_{x<M}\int_{0}^{x}d\pi\int_{x}^{M}\frac{1}{\nu(t)}dt.$$

If $B^{-}_{M}$ et $B^{+}_{M}$ are bounded,  then the following Poincar\'e inequality holds 
$$\mbox{Var}_\pi(f)\le c_P\int f'^2d\nu.$$

In addition, the optimal constant $c_P$ verifies $$\frac{1}{2}\max(B^{+}_{M},B^{-}_{M})\leqslant c_P\leqslant 4\max(B^{+}_{M},B^{-}_{M})$$.

\end{lemma}

We apply the lemma to $\pi$ and $\nu=x(1-x)\pi$ using both sides of the estimates. Denote $\pi_0$ and $\nu_0$ the case where $s=0$ and the Poincar\'e constant is $m$. Then
\begin{eqnarray*}
B^{+,s}_M&=&\sup_{x>M}\int_x^1d\pi\int_M^x\frac{1}{\nu(t)}dt\\
&=&\sup_{x>M}\int_x^1 e^{st}d\pi_0(t)\int_M^x\frac{e^{-st}}{\nu_0(t)}dt\\
&\le& e^{(1-M)s}B^{+,0}_M
\end{eqnarray*}
The same reasoning shows that $B^{-,s}_M\le B^{-,0}_M$.
\end{proof}
Of course, one can do easily the same for $s<0$ using a symmetric reasoning. If the order is good with respect to the immigration parameter, as the case $s=0$ is optimal, it is an open question to look at the dependence with respect to the selection parameter. We may also consider a convergence in entropy, via the logarithmic Sobolev inequality without selection established by Stannat \cite{stannat2000} or Miclo \cite{miclo2003} and the same line of proof using Holley-Stroock perturbation argument or the Hardy type condition for logarithmic Sobolev inequality (see again \cite{ABC}). Note that the  convergence in entropy entails a convergence in total variation via Cszisar-Pinsker-Kullback inequality, but the constant involved are less explicit so we omit the details.\\\medskip

This quantitative long time behaviour enables us to give an error while approximating the asymptotic Simpson index (being a smooth function of the species). As usual, an $L^2$ decay will enable us to consider long time behaviour for initial measures whose density with respect to the invariant measure is bounded, which in could prevent starting from a Dirac measure. However due to regularization, and so waiting a time $t_0$, enables (loosing on the constants in the decay) to start from a Dirac measure. See for example \cite{BCG}.

\section{Generalization to a larger number of species or in random environment}
In this section we provide extensions of the two species case to 1) finite number of species, 2) two species case in a particular random environment, namely Wright-Fisher diffusion environment.

\subsection{Expectation approximation for three species}
\quad\\
In fact we will give the main ideas for $S=2$. Extension to a larger number of species is only technically involved and requires no further arguments. Denote $X_{t}$ and $Y_{t}$ the proportions of the two main species, $s^{x}_{t}$ and $s_{t}^{y}$ their selection parameters and $p_{t}^{x}$ and $p_{t}^{y}$  their proportions in the pool. The immigration parameter will still be denoted $m_t$. The method presented for $ S = 1 $ in the previous section can be generalized to a larger number of species. It will have of course some limitations: greater the number of species is, larger will be the size of the approaching linear system. In fact, the derivative of the expectation of order $ n $ involves only the expectation of lower and higher order in the case of two species. Now with 3 species we also need to know the expectation of the form $ \E [X_{t}^{n}Y_{t}^{k}] $ for $ k, n$ in $\mathbb{N} $. We will thus need a system of size $ N^{2} $.
 We present here the extension of our approximation for 3 species.

\begin{equation}
\begin{pmatrix}
dX_{t}\\
dY_{t}
\end{pmatrix}
=
\begin{pmatrix}
m_{t}(p_{t}^{x}-X_{t})+X_{t}(s_{t}^{x}-X_{t}s_{t}^{x}-Y_{t}s_{t}^{y})\\
m_{t}(p_{t}^{y}-Y_{t})+Y_{t}(s_{t}^{y}-X_{t}s_{t}^{x}-Y_{t}s_{t}^{y})\\
\end{pmatrix}dt
+ \sigma (X_{t},Y_{t})dB_{t} \label{eqprincipalemulti}
\end{equation}

where $\sigma$ verifies $\sigma.\sigma^{*}(x,y)=a(x,y)$ with \\
$$a(x,y)=2\begin{pmatrix} x(1-x)&-xy\\ -xy&y(1-y)\\ \end{pmatrix} $$

We have to calculate $d(X_{t} ^{n}Y_{t}^{k})$  with It\^o's formula :

\begin{align}
d(X_{t} ^{n}Y_{t}^{k})= \nonumber &\big(m_{t}(p_{t}^{x}-X_{t})+X_{t}(s_{t}{x}-X_{t}s_{t}^{x}-Y_{t}s_{t}^{y})\big)nX_{t} ^{n-1}Y_{t}^{k}\\ \nonumber
&+\big(m_{t}(p^{y}_{t}-Y_{t})+Y_{t}(s^{y}_{t}-X_{t}s_{t}^{x}-Y_{t}s^{y}_{t})\big)kX_{t} ^{n}Y_{t}^{k-1}\\ \nonumber
&+n(n-1)(1-X_{t})X_{t} ^{n-1}Y_{t}^{k} +k(k-1)(1-Y_{t})X_{t} ^{n}Y_{t}^{k-1}\\ \nonumber
&-2nkX_{t} ^{n}Y_{t}^{k}+ d\mathscr{M}_{t}\\ \nonumber
=&X_{t} ^{n-1}Y_{t}^{k}n(m_{t}p_{t}^{x}+n-1)\\ \nonumber
&+X_{t} ^{n}Y_{t}^{k-1}k(m_{t}p^{y}_{t}+k-1)\\ \nonumber
&+X_{t} ^{n}Y_{t}^{k}\big(-m_{t}(n+k)-2kn-k(k-1)-n(n-1)+ns_{t}^{x}+ks_{t}^{y}\big)\\ \nonumber 
&-X_{t} ^{n+1}Y_{t}^{k}s_{t}^{x}(n+k)\\ \nonumber
&-X_{t} ^{n}Y_{t}^{k+1}s_{t}^{y}(n+k)\\
&+d\mathscr{M}_{t}\label{coeff3speces}
\end{align}
here $\mathscr{M}_{t}$ is a martingale. Then $ d\E [X_ {t} ^ {n} Y_ {t} ^ {k}] $ is expressed in terms of 4 other quantities which complicates the one dimensional calculations. Moreover we must define what are the neglected expectations on which we will make an approximation, that is how to close the system. We can decide  we make an approximation to the order $ N $ then that we neglect all the terms of higher order in the expression of $dE[X_{t} ^{n}Y_{t}^{k}]$ where $max(n,k)=N$.

Suppose we want to get the expectation up to order $N$ we need exactly $\sum_{k=1}^{N}2k+1$ quantities. And so the size of the approaching differential linear system will be of order $ (N + 1) ^ {2} $. \\
The following figure represents the complexity of the problem, for example $ d\E [X_{t}Y_{t}] $ is expressed as a function of the expectations of the quantities to which the blue arrows point.
$$\begin{matrix}
\quad X\textcolor{orange}{\rightarrow} \qquad \textcolor{blue}{\leftarrow}XY \textcolor{blue}{\rightarrow} \qquad Y\\
\quad \textcolor{orange}{\swarrow} \quad   \qquad  \quad \textcolor{blue}{\swarrow}\quad \quad \quad \textcolor{blue}{\searrow}\qquad \qquad\\
\qquad X^2\qquad X^{2}Y \qquad X^{2}Y^{2}\qquad Y^{2}X\qquad Y^{2}
\end{matrix}$$

Algorithmically it is not very difficult to build the matrix approaching the expectations of the diffusion. We must begin by giving a vector composed of the different expectation of size $ (N + 1)^2 $. For that, let us define an application that transforms the expectation of order $ (n, k) $ that is to say $ \E [X_{t} ^ {n} Y_{t} ^{k}] $ into an integer which corresponds to its coordinate in the expectation vector. 
Next we build the matrix $A_{N²}$ as in the case of two species from the coefficients calculated in \ref{coeff3speces}.
We consider as error each expectation $(n,k)$ with $max(n,k)>N$ in the It\^o formula. So that the error is composed of  $2N+1$ terms. And the approximation boils down to solving numerically a linear system. As an example consider the case  $N=1$, we therefore involve four expectations which are $\E[X^{n}Y^{k}], k,n \in \{0,1\} $. The order imposed by $ \phi $ is therefore $(1,\E[X],\E[Y],\E[XY])$. Three terms compose the error: $\E[s_{t}^{x}X^{2}],\E[s_{t}^{y}Y^{2}],2\E[(s_{t}^{x}+s_{t}^{y})XY^{2}]$. We can also, as in the case of two species, prove the convergence of this algorithm by following exactly the same pattern as in the one species case.
The renormalizing coefficients of the expectation $(n,k)$ then become
$\frac{(s_{t}^{x})^{n}(s_{t}^{y})^{k}}{(n-1)!(k-1)!}$. It can similarly be shown that the error is at most of the order of $\frac{N^{2}\max(s_{t}^{x},s_{t}^{y})^{N}}{(N-1)!}$.\\

The extension to a larger number is straightforward and will entail an error of the order $\frac{N^{\frac{S}{2}+2}\|s_t\|_\infty^N}{N!}$, and it will still be reasonable but requires computations of a system of size $N^{S}$ which may be prohibitive for large $S$.

\subsubsection*{Numerical applications}

We can easily program such an algorithm and check that the results obtained are in agreement with quantities obtained by Monte Carlo method. See following figures:

\underline{Basic example}

\begin{figure}[htbp]
\label{basic}
\begin{minipage}[c]{.45\linewidth}
\begin{center}
\includegraphics[scale=0.35]{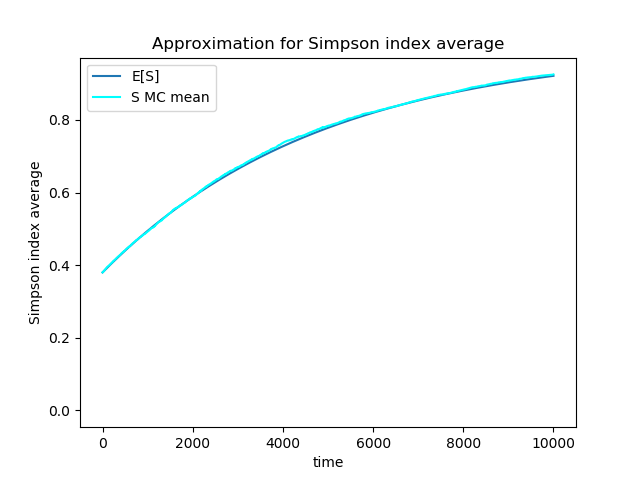}
\end{center}
\end{minipage}
\hfill
\begin{minipage}[c]{.45\linewidth}
\begin{center}
\includegraphics[scale=0.35]{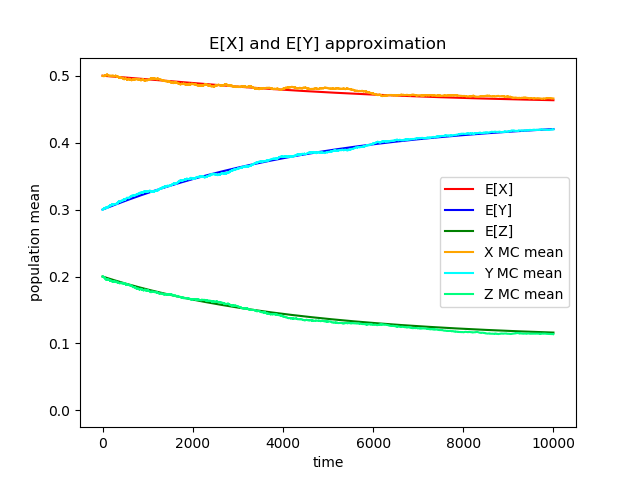}
\end{center}
\end{minipage}
\label{fig:image7}
\caption{On the left hand side expectation of the Simpson index and Monte Carlo mean, and on the right hand side expectations of each species and their Monte Carlo means}
\end{figure}
We consider here a case with no immigration and constant selection parameter. The number of simulated trajectories for  MC mean is  1000, $J=1000$, $m=0$, $X_{0}=0.5$, $Y_{0}=0.3$, $s_{y}=2$, $s_{x}=1$, the size of the approaching linear system is 144. Figure 8 plots approximate values of $\E[\mathscr{S}_{t}]$ and $\E[X_{t}]$ by the precedent method from the approximation in large population and by MC method from the discrete model.
\quad\\

\underline{Time dependent parameter case}

\begin{figure}[htbp]
\begin{minipage}[c]{.45\linewidth}
\begin{center}
\includegraphics[scale=0.35]{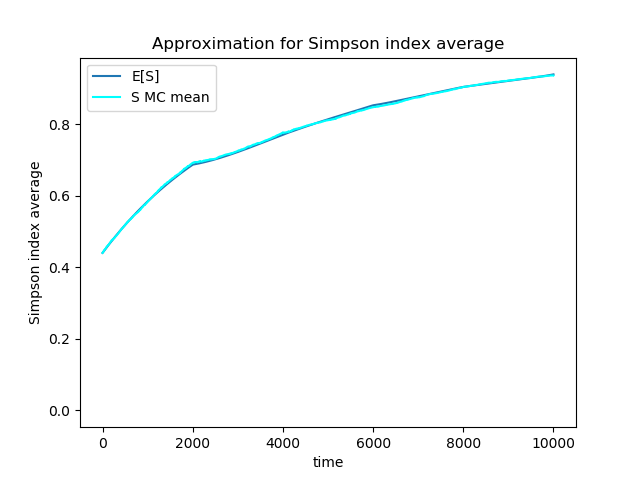}
\end{center}
\end{minipage}
\hfill
\begin{minipage}[c]{.45\linewidth}
\begin{center}
\includegraphics[scale=0.35]{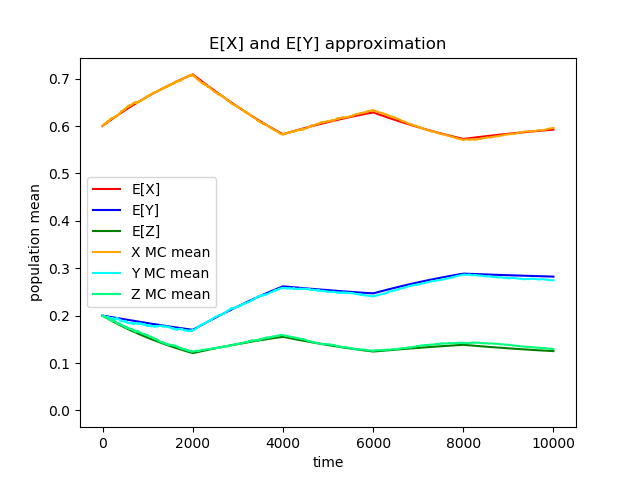}
\end{center}
\end{minipage}
\label{fig:image11}
\caption{Left hand side Simpson index, right hand with the expectations of the three species.}
\end{figure}
In this second example, we consider once again a case without immigration and time dependent selection parameter. The number of simulated trajectories for  MC mean is  1000, $J=1000$, $m=0$, $X_{0}=0.5$, $Y_{0}=0.3$, $s_{y}=2$, $s_{x}$ is piecewise constant taking two values $4$ and $-4$ at regular time intervals, the size of the approaching linear system is 144.  Figure 9 plot the approximate values of $E[\mathscr{S}_{t}]$ and $E[X_{t}]$ by the precedent method from the approximation in large population and by MC method from the discrete model.

\subsection{When the selection is a diffusion.}

In the third Section we gave a method to get the moments of  $X_{t}$, and thus  $\E(\mathscr{S}_{t})$ for a time dependent immigration/selection parameter. If these parameters are random but autonomous, it gives a way to approximate the expectation of the Simpson index by doing a Monte Carlo mean with respect to the environment, passing from quenched to annealed. It would be however more interesting to evaluate directly the expectation of the Simpson index without further Monte Carlo simulations. It seems quite impossible to give a general algorithm for every environment but we will give in this section an efficient approximation method in a particular case. We consider for the selection parameter $s_t$ a rescaled Wright-Fisher diffusion, whose leading Brownian motion is independent of the one leading the SDE for the species evolution. This choice assures us that  $s_t$  is a diffusion evolving in a bounded set and the choice of the different parameters leads to a wide choice of a Moran process with immigration.

\subsubsection{The expectation approximation}
\quad\\
Let us just give first the diffusion approximation result for this particular case, whose proof is even simpler as it relies on usual approximation diffusion for Markov chains.
\begin{theorem}
Assume that $(v^{J}_{n})_{n \in \mathbb{N}}$ is a Moran process without selection with size $J$ and the parameters $m^{s}$ et $p^{s}$.
Let $c$ and $b$ two constants such as  $s^{J}_{n}= cv^{J}_{n}-b \quad  \forall n \in \mathbb{N}$, and assume that  $X_{n}$ follow a Moran process with size $J$ and parameters $(m_{n})_{n \in \mathbb{N}}$,  $(p_{n})_{n \in \mathbb{N}}$ et  $(s^{J}_{n})_{n \in \mathbb{N}}$ describe in the first part.
Let $U^{J}_{n}$ be the process having for coordinates  $X^{J}_{n}$ et $s^{J}_{n}$.

Then when $J$ goes to infinity, the process $U^{J}_{tJ²}$ converge in law to the process $U_{t}$ which coordinates are solutions of the following stochastic differential equation:
\begin{equation}
\begin{pmatrix}
dX_{t}\\
dv_{t}
\end{pmatrix}
=
\begin{pmatrix}
m'_{t}(p_{t}-X_{t})+s'_{t}X_{t}(1-X_{t})\\
m'^{s}_{t} (p^{s}_{t}-v_{t})
\end{pmatrix}dt
+\begin{pmatrix}
\sqrt[]{2X_{t}(1-X_{t})}\\
\sqrt[]{2v_{t}(1-v_{t})}
\end{pmatrix} dB_{t} \label{eqprincipale} 
\end{equation}

where $s_{t}'=\frac{s_{t}}{J}=cv_{t}-b$,  $m_{t}'=\frac{m_{t}}{J}$, $m_{t}'^{s}=\frac{m_{t}^{s}}{J}$.
\end{theorem}

To approach the expectation of $X_{t}$ we use the method describe previously for three species, here  $v$ play the same role as a third species.
However the dynamics is not exactly the same, the It\^o formula gives us:

\begin{align}
d(X_{t} ^{n}v_{t}^{k})= \nonumber &X_{t} ^{n-1}v_{t}^{k}n(m'_{t}p_{t}+n-1)\\ \nonumber
&+X_{t} ^{n}v_{t}^{k-1}k(m'^{s}_{t}p^{s}_{t}+k-1)\\ \nonumber
&+X_{t} ^{n}v_{t}^{k}\big(-(m'_{t}+b)n-km_{t}'^{s}-k(k-1)-n(n-1)\big)\\ \nonumber 
&+X_{t} ^{n+1}v_{t}^{k}nb\\ \nonumber
&+X_{t} ^{n}v_{t}^{k+1}cn\\
&-X_{t} ^{n+1}v_{t}^{k+1}nc\\
&+d\mathscr{M}_{t}\label{coeff3speces2}
\end{align}
with $\mathscr{M}_{t}$ a martingale. Then as previously we close our system, for a given $N$, and to do so to neglect all the terms of higher order in the expression of $d\E[X_{t} ^{n}v_{t}^{k}]$ where $max(n,k)=N$. And now the algorithm is able to calculate all the expectations of the form $ \E[X_{t} ^{n}v_{t}^{k}]$ and so obtain the expectation of the Simpson index. The proof follows the same pattern. The renormalizing coefficients of the expectation $(n,k)$ in the proof allow to control the eigenvalues of the matrix thanks to the Gershgorin disks as before.
Many choices  are possible and we take here the coefficient $\frac{1}{\sqrt[4]{n!k!}}$. This choice  leads to a convergence speed at most  of the order of $\frac{N^2}{\sqrt[4]{(N)!}}$.

\begin{figure}[htbp]
\begin{minipage}[c]{.45\linewidth}
\begin{center}
\includegraphics[scale=0.35]{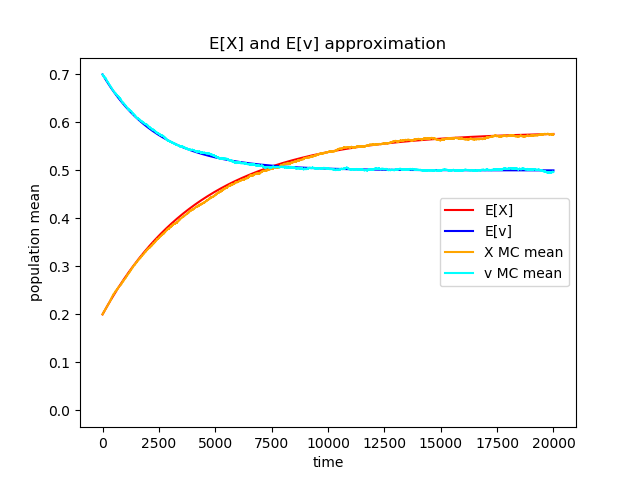}
\end{center}
\end{minipage}
\label{fig:image12}
\caption{Are plotted the approximate values of $\E[v_{t}]$ and $\E[X_{t}]$ by the precedent method from the approximation in large population and by MC method from the discrete model. The number of simulated trajectories for  MC mean is  $5000 $, $J=1000$, $X_{0}=0.2$, $v_{0}=0.7$, $m^{s}=4$, $m=2$, $p=p^{s}=0.5$, $c=3$, $b=0.5$ the size of the approaching linear system is 144}
\end{figure}

\subsubsection{Comparison with the neutral model.}

\quad\\
In this part we compare the case where $s$ is "neutral on average",  to the neutral case with $s=0$.
Thanks to the previous method one can for example calculate the average Simpson index in the case where the selection expectation is 0. For it let's take $p^{s}=1/2$, $v_{0}=1/2$ (this enforces $c=-2b$). The following figures show the results:
\begin{figure}[htbp]
\begin{minipage}[c]{.45\linewidth}
\begin{center}
\includegraphics[scale=0.3]{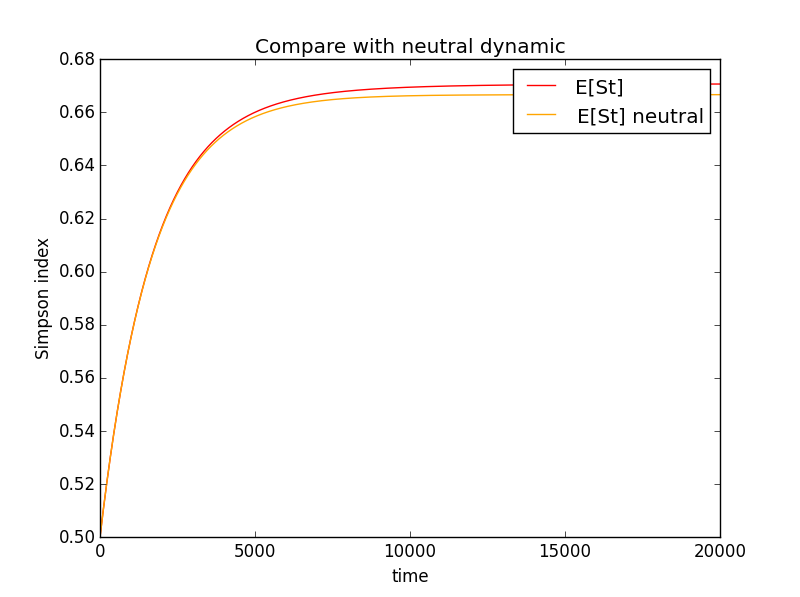}
\end{center}
\end{minipage}
\hfill
\begin{minipage}[c]{.45\linewidth}
\begin{center}
\includegraphics[scale=0.3]{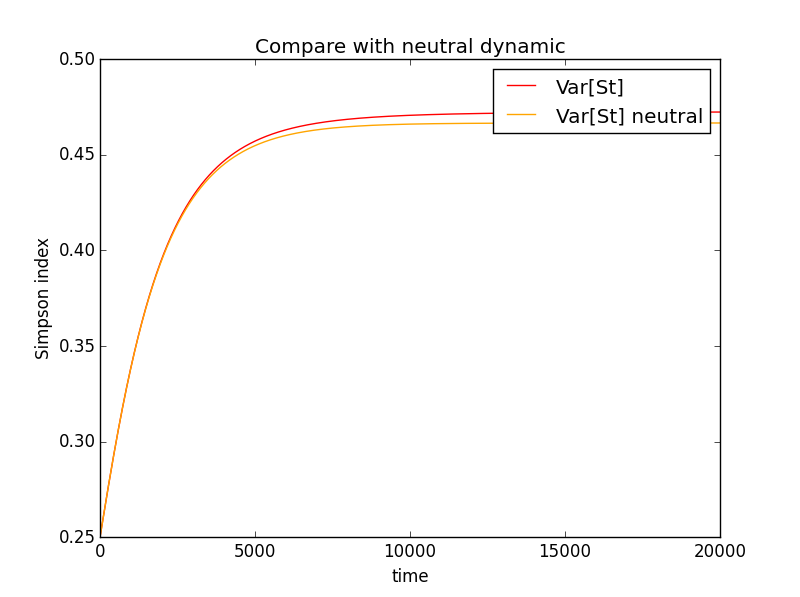}
\end{center}
\end{minipage}
\label{fig:image13}
\caption{Comparison with neutral case. Approximate values of $\E[\mathscr{S}_{t}]$ and $\E[X_{t}]$ by the precedent method from the approximation in large population  $X_{0}=0.5$, $m^{s}=1$, $m=2$, $p=0.5$, $c=3$, size of the approaching linear system is 144}
\end{figure}

\begin{figure}[htbp]
\begin{minipage}[c]{.45\linewidth}
\begin{center}
\includegraphics[scale=0.3]{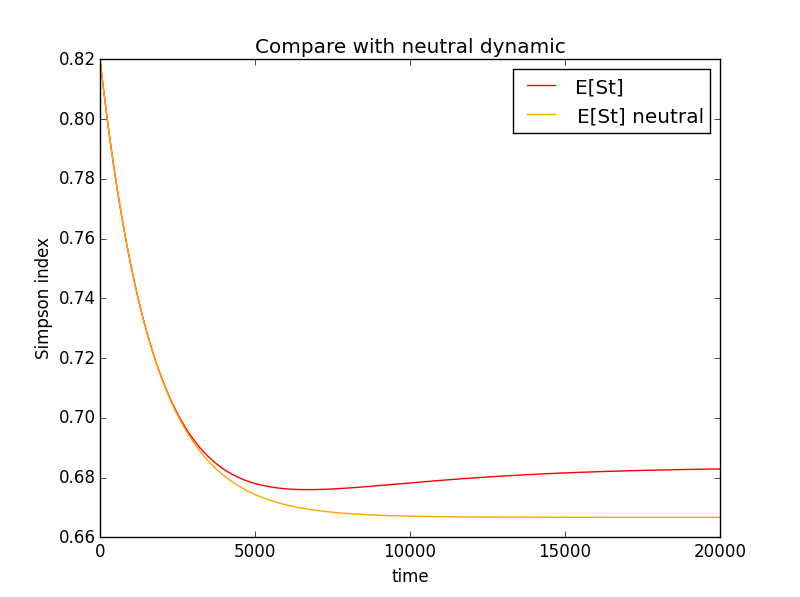}
\end{center}
\end{minipage}
\hfill
\begin{minipage}[c]{.45\linewidth}
\begin{center}
\includegraphics[scale=0.3]{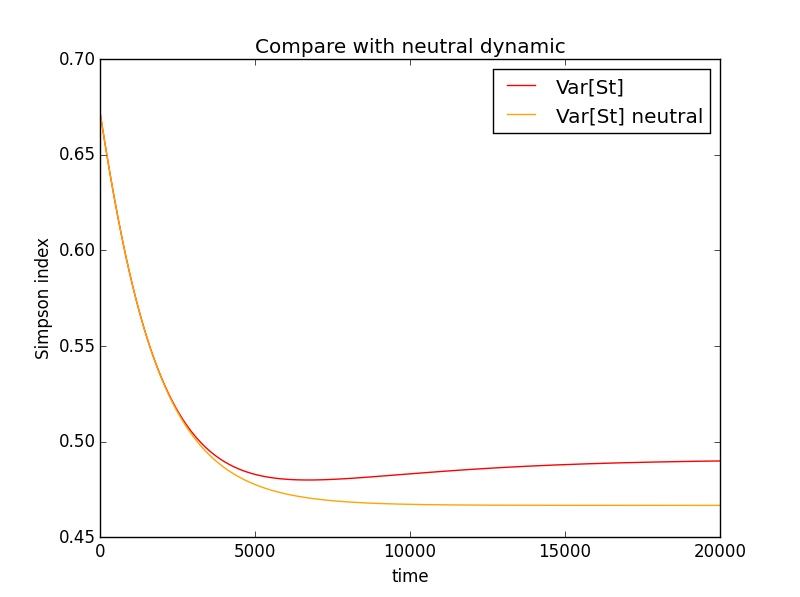}
\end{center}
\end{minipage}
\label{fig:image14}
\caption{Comparison with the neutral case. Approximate values of $\E[\mathscr{S}_{t}]$ and $\E[X_{t}]$ by the precedent method from the approximation in large population  $X_{0}=0.1$, $m^{s}=0.5$, $m=2$, $p=0.5$, $c=5$, size of the approaching linear system is 144.}
\end{figure}

We thus see that a selection even if neutral in mean, involves deeper mechanism which  lead to a different behaviour than the neutral one. Of course the Simpson index involves not only the expectation of one species but also the moment of order two.

\subsection{Effect of selection on increase of biodiversity}
\quad\\
We have already seen in the case of two species that selection alone could contribute to the decrease of the average Simpson index in the absence of immigration. There was however a threshold for $ s $ under which such a phenomenon could not occur. We sort of generalize it here to any number of species.

\begin{prop}
Note as previously $ S + 1 $ the number of species in the community with $s_i$ the selection parameter for species $i$. Then if all $ \|s_ {i} \|_{\infty} $ are less than $ \frac {1}{2} $ the Simpson index is increasing in the absence of immigration.
In other words, the selection can not be the source of the diversity decreasing.
\end{prop}

\begin{proof}


Assume all  the  $s_{t}^{i}, \forall t$ are between $\alpha$ and $-\alpha$ .
 First write:
\begin{align*}
1-\mathscr{S}_{t}=2&\sum\limits_{i=1}^{S}X^{i}_{t}(1-X^{i}_{t})-\sum\limits_{j\neq i}X^{i}_{t}X^{j}_{t}\\
=&\sum\limits_{i=1}^{S}X^{i}_{t}(1-X^{i}_{t})+\sum\limits_{i=1}^{S}X^{i}_{t}(1-\sum\limits_{i=1}^{S}X^{i}_{t})\\
=&\sum\limits_{i=1}^{S+1}X^{i}_{t}(1-X^{i}_{t})\\
\end{align*}

Then, 
\begin{align*}
d\E[\mathscr{S}_{t}]=&2\E[1-\mathscr{S}_{t}]-2\E[\sum\limits_{i=1}^{S}s_{t}^{i}X_{t}^{i}(\mathscr{S}_{t}-X_{t}^{i})]dt\\
= & 2\E[1-\mathscr{S}_{t}-\sum\limits_{i=1}^{S}s_{t}^{i}X_{t}^{i}(\mathscr{S}_{t}-1)-\sum\limits_{i=1}^{S}s_{t}^{i}X_{t}^{i}(1-X_{t}^{i})]dt\\
\geqslant & 2\E[1-\mathscr{S}_{t}-\sum\limits_{i=1}^{S}s_{t}^{i}X_{t}^{i}(\mathscr{S}_{t}-1)-\alpha\sum\limits_{i=1}^{S}X_{t}^{i}(1-X_{t}^{i})]dt\\
\geqslant & 2\E[(1-\mathscr{S}_{t})(1+\sum\limits_{i=1}^{S}s_{t}^{i}X_{t}^{i})-\alpha\sum\limits_{i=1}^{S}X_{t}^{i}(1-X_{t}^{i})]dt\\
\geqslant & 2\E[(1-\mathscr{S}_{t})(1-\alpha+\sum\limits_{i=1}^{S}s_{t}^{i}X_{t}^{i})+\alpha\sum\limits_{i=1}^{S}X_{t}^{i}(1-\sum\limits_{i=1}^{S}X_{t}^{i})]dt\\
\geqslant & 2\E[(1-\mathscr{S}_{t})(1-\alpha+\sum\limits_{i=1}^{S}s_{t}^{i}X_{t}^{i})]dt\\
\geqslant & 2\E[(1-\mathscr{S}_{t})(1-\alpha(1+\sum\limits_{i=1}^{S}X_{t}^{i}))]dt
\end{align*}

and so if $\forall i,  \alpha\leqslant \frac{1}{2}$ , $ d\E[\mathscr{S}_{t}]\geqslant 0$ and $\E[\mathscr{S}_{t}]$ is increasing.

\end{proof}
Remark that this bound is certainly not optimal, as the two species case indicates but true for each $S$.

%
%

\subsection{Long time behaviour.}
\quad\\
We will once again assume in this part $s,m,p$  are constants. If $m=0$ then a species will still  invade the community definitively. On the other hand, if $ m \neq 0 $, the law of the vector of abundance, converges in a long time to a unique invariant measure. Consider the generator of  the diffusion \eqref{eqprincipale} which is the generator of the Wright-Fisher diffusion with selection and mutation:
$$\mathscr{L}f(x)=\sum\limits_{i,j=1}^{S}x^{i}(\delta_{i,j}-x^{j})\frac{\partial²f(x)}{\partial x^{i}\partial x^{j}j}+\sum\limits_{i=1}^{S}\left(m(p^{i}-x^{i})+x^{i}\left(s^{i}-\sum\limits_{i=1}^{S}x^{i}s^{i}\right)\right)\frac{\partial f(x)}{\partial x^{i}}$$
A reversible and stationary measure for the diffusion \eqref{eqprincipale} is given by (see for example \cite{theinfinitneutralalleles,wright1977,Griffiths}:

$$\pi_{S}(dx)=C\times \exp\left(\sum\limits_{i,j=1}^{S+1}s^{i}x^{i}x^{j}\right)\times (x^{1})^{mp^{1}-1}\times...\times (x^{S+1})^{mp^{S+1}-1}dx_{1}...dx_{S}$$

Where  $x^{S+1}=1-\sum\limits_{i=1}^{S}x^{i}$ and $p^{S+1}=1-\sum\limits_{i=1}^{S}p^{i}$, $s^{S+1}=0$. C is a constant just like $\int \pi_{S}(dx)=1$. 

Of course, when $s$ and $m$ are  time dependent, periodic for example, an invariant measure will not exist. The next figure presents the approximate values of $\E[\mathscr{S}_{t}]$ and $\E[X^{i}_{t}]$ for $i$ in $\{1,2,3\}$ by the precedent method from the approximation in large population and by Monte Carlo method from thr discret model. The number of simulated trajectories for  Monte Carlo mean is  $5000$, $J=500$, $m$ is a time dependant piecewise process, it takes alternatively the values of $3$ and $0$ at regular time intervals. $X_{0}=0.5$, $p_{x}=0.33$, $Y_{0}=0.3$, $p_{y}=0.33$, $s_{y}$ et  $s_{x}$ are Markovian jump processes, the size of the approaching linear system is 144.

\begin{figure}[htbp]
\begin{minipage}[c]{.45\linewidth}
\begin{center}
\includegraphics[scale=0.35]{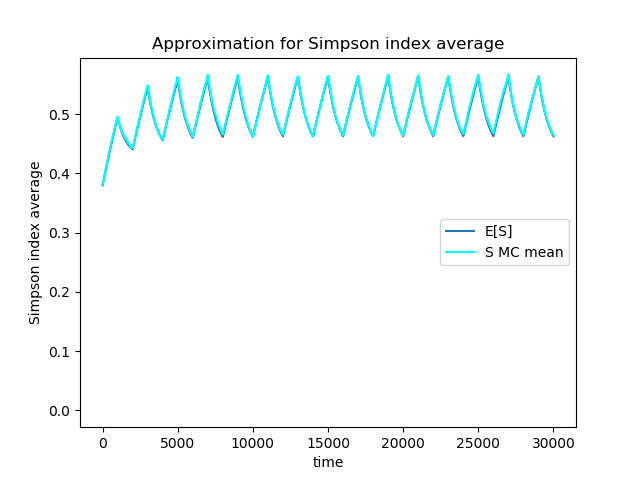}
\label{malea}
\end{center}
\caption{$\E[\mathscr{S}_{t}]$}
\end{minipage}
\hfill
\begin{minipage}[c]{.45\linewidth}
\begin{center}
\includegraphics[scale=0.35]{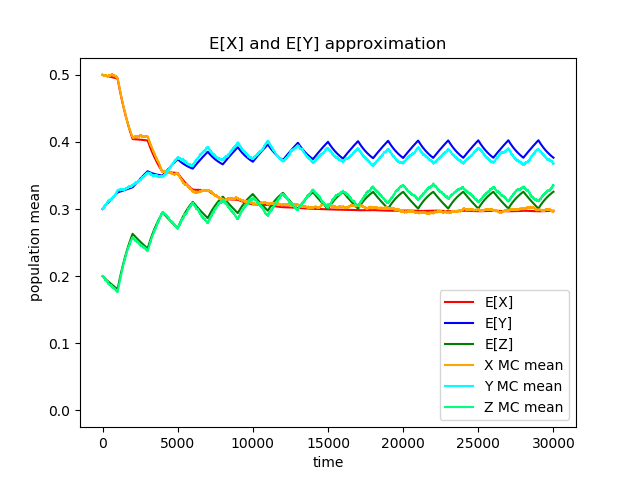}
\label{}
\end{center}
\caption{$\E[X_{t}], \E[Y_{t}]$}
\label{Tux}
\end{minipage}
\end{figure}

Concerning the long time behaviour, we may once again refer to \cite{shimakura1977} for the spectral gap which is $e^{(S+1)\sum_1^{S+1}s^i}/m$ by Holley-Stroock's perturbation argument. Unfortunately, it is not possible to refine this argument as there is no Hardy's type inequalities in this case. Once again it is also possible to derive a logarithmic Sobolev inequality, and thus convergence in entropy (and total variation) but constants are less explicit.

\section{Proofs}
In this section we gather the proofs, technical or more or less well known.

\subsection{Proof of the diffusion approximation, Theorem \ref{th:diffapprox}}
\label{diffapprox}
\quad\\

In the following proof we'll get back to a  martingale problems. All the results used in this section can be funded in \cite{diffmultidim} p267-272.

For the sake of clarity, assume that $ m = 0 $, and that $ S = 1 $. \\The multidimensional case is treated exactly the same way. \\
We can put  $h=\frac{1}{J²}$ and $U$ means here $(x,s)$ where $x \in E_{x}$ and $s \in E_{s}$.\\
Let $a(U)=x(1-x)$,  $b(U)=sx(1-x)$ and $L_{x}f(U)=b(U)\times \frac{\partial f}{\partial x}(U)+a(U)\times \frac{\partial² f}{\partial x²}(U)$ the generator of  the SDFE \eqref{eqprincipale} and $L_{s}f(U)=\sum\limits_{y\in E_{s}}
Q_{s,y}f(x,y)$ the generator of a Markovian jump process applied to a function depending of the population variable.

Let's start with the following lemma
\begin{lemma}(\cite{diffmultidim} p268)
\quad\\
Let $f$ be a  $C^{\infty}$ function, note $A_{J}f(U)=\int_{E}f(y)-f(u)d\pi_{J}(U,dy)$ then $J²A_{J}f$ converge uniformly to $L_{x}f+L_{s}f$
\end{lemma}
\begin{proof}
\begin{align*}
A_{J}f(U)=&\int_{E}f(y)-f(u)d\pi_{J}(U,dy)=\int_{E}f(z,w)-f(x,s)d\pi_{J}(U=(x,s),dy)\\
=&\int_{E}f(z,w)-f(x,s)d\pi_{J}(U,dy)\\=&\int_{E}f(z,w)-f(x,w)+f(x,w)-f(x,s)d\pi_{J}(U,dy)\\
=&\int_{E}f(z,w)-f(x,w)d\pi_{J}(U,dy)+\int_{E}f(x,w)-f(x,s)d\pi_{J}(U,dy)\\
=&\int_{E}  f(z,w)-f(x,w)d\pi_{J}(U,dy)		+	\sum\limits_{w\in E_{s}}f(x,w)-f(x,s)P_{J}^{i}(x,y,h).
\end{align*}
 
Via Taylor's formula, we obtain
\begin{align*}
&\int_{E}  f(z,w)-f(x,w)d\pi_{J}(U,dy)=\E[f(X_{t+h},s_{t+h})-f(x,s_{t+h})|U_{t}=(x,s)]\\
&=\frac{\partial f}{\partial x}(x,s)\E[X_{t+h}-x|U_{t}]+\frac{\partial² f}{2\partial x²}(x,s)\E[(X_{t+h}-x)²|U_{t}]\\ 
&+ \frac{\partial² f}{\partial x\partial s}(x,s)\E[X_{t+h}-x|U_{t}]E[s_{t+h}-s|U_{t}]+O(\E[\| U_{t+h}-U_{t}\|³|U_{t}]).
\end{align*}

we give the limits of the previous quantities,
\begin{itemize}
\item $\lim\limits_{h\rightarrow 0} \sup\limits_{E}\|\frac{1}{h}\E[X^{J}_{t+h}-x|U^{J}_{t}]\|=b(x,s)$ by the second propriety (\ref{exp-var}),
\item $\lim\limits_{h\rightarrow 0} \sup\limits_{E}\|\frac{1}{h}\E[(X^{J}_{t+h}-x)²|U^{J}_{t}]\|=a(x,s)$ by the second propriety (\ref{exp-var}),
\item $\lim\limits_{h\rightarrow 0} \sup\limits_{E}\|\frac{1}{h}\E[X^{J}_{t+h}-x|U^{J}_{t}]\E[s^{J}_{t+h}-s|U^{J}_{t}]\|=0$  because
\begin{eqnarray*}
|\E[s^{J}_{t+h}-s|U_{t}]|&=&|\sum\limits_{w\in E_{s}}(w-s_{t})P_{J}^{i}(s_{t},w,h)|\\& \leqslant& \max\limits_{x,y\in E_{s}}|x-y|\sum\limits_{w\in E_{s}/s_{t}}P_{J}^{i}(s_{t},w,h),
\end{eqnarray*}
\item $\lim\limits_{h\rightarrow 0} \sup\limits_{E}\E[\| U^{J}_{t+h}-U\|³|U^{J}_{t}]=0$  because $\int_{E}\sup\limits_{i\leqslant S}|y_{i}-u_{i}|^{3}\pi_{J}(U,dy) \rightarrow 0 $  and by the previous point.
\end{itemize}

Then, going to the limit in the previous expression, 
\begin{align*}
\lim\limits_{h\rightarrow 0} &\sup\limits_{E}|A_{J}f(U)\frac{1}{h}-L_{x}f(U)-L_{s}f(U)|\\
\le&\lim\limits_{h\rightarrow 0} \sup\limits_{E}|\frac{1}{h}\int_{E}  f(z,w)-f(x,w)d\pi_{J}(U,dy)-	L_{x}f(U)|\\ 
&\qquad+\lim\limits_{h\rightarrow 0} \sup\limits_{E}|\frac{1}{h}\sum\limits_{w\in E_{s}}\left(f(x,w)-f(x,s)\right)P_{J}^{i}(x,y,h)-L_{s}f(U)|\\ =&0
\end{align*}
And this expression conclude the proof.

\end{proof}
Now let $f$ be $\mathscr{C}^{\infty}$, then 
\begin{align*}
\E_{U}[f(U^{J}_{t})]&=f(U)+\E_{U}\left[  \sum\limits_{k=1}^{|tJ²|-1}\E[f(U^{J}_{(k+1)h})-f(U^{J}_{kh})|U^{J}_{kh}]\right]\\
&=f(U)+\E_{U}\left[  \sum\limits_{k=1}^{|tJ²|-1}A_{J}f(U_{kh})\right]\\
&=f(U)+\E_{U}\left[  \sum\limits_{k=1}^{|tJ²|-1}\frac{h}{h}A_{J}f(U_{kh})\right]\\
\end{align*}
and so  $$\E_{U}[f(U^{J}_{t+h})-f(U)-\sum\limits_{k=1}^{|tJ²|-1}\frac{h}{h}A_{J}f(U_{kh})]=0$$
i.e $f(U^{J}_{t})-f(U)-\sum\limits_{k=1}^{|tJ²|-1}\frac{h}{h}A_{J}f(U_{kh})$ is a  martingale for $ \pi^{J}$.\\
Moreover, note the sum is a Riemann sum and the previous lemma  ensures when $ J $ tends to infinity the convergence of $$f(U^{J}_{t})-f(U)-\sum\limits_{k=1}^{|tJ²|-1}\frac{h}{h}A_{J}f(U_{kh})$$ towards $$f(U_{t})-f(U)-\int_{0}^{t}L_{x}f(U_{s})+L_{s}f(U_{s})ds.$$ 
We need now to find a probability measure on the Borel sets of the canonical space  $\mathscr{C}([0,1],R)$ verifying the martingale problem for $ L_ {x} + L_ {s} $.
Let us show now that $ \pi ^ {J} $ admits an adherent value in the space of probability measure on the Borel of $\mathscr{C}([0,1],R)$ with the norm $$\|\pi^{J}\|=\sup\limits_{f\in C
}\frac{|\int fd\pi^{J}|}{|f|_{\infty}}$$ (which is a norm since the $ \pi^{J} $ are supported in $ [0,1] $).\\
Let note $ \pi^{J} f=\int_{E} fd\pi^{J}$.

Let $ (f_{n}) _ {\mathbb{N}} $ be a dense sequence in the space of continuous functions then $ (\pi^{J} f_ {n}) _ {J} $ is a sequence of $ \mathbb{R} $ having an adherence value in $ \mathbb {R} $ because it is uniformly bounded by $ |f_{n}|_{\infty} $. Then by diagonal extraction, eventually for a subsequence, $ \pi ^ {J}f_{n} $ converges to a certain $ \phi_ {fn} $ in $ \mathbb{R} $ for all $ n $. And by the uniformly continuous extension theorem, we define $ \phi_ {f} $ for all $ f $ of $ \mathscr{C} $. And since $ \phi_ {f} $ is a linear form, the Riesz-Markov theorem ensures the existence of a unique measure $ \mu $ such that $ \int fd \mu = \phi_ {f} $.
Since this is true for all $ f \in \mathscr{C} $, by considering the constant function equal to $ 1 $, we find $ \mu (\Omega) = 1 $ and $ \mu $ is a probability. The convergence of $ \| \pi^{J} - \mu \| $ to $ 0 $ is then immediate in view of the chosen norm.\\
Thus our sequence $ \pi^{J} $ admits an accumulation point. So, there is at least one $ \pi $ and one $ X $ process that satisfy the martingale problem associated with $ L_ {x} + L_ {s} $. \\
And so  $\pi$ verifies $\int f(U_{t})-f(U)-\int_{0}^{t}L_{x}f(U_{s})+L_{s}f(U_{s})ds d\pi(U,U_{t})=0$. 
So it exist at least a solution to the martingale problem for $ L_ {x} + L_ {s} $ .

\quad\\ 
If the uniqueness of this martingale problem is verified then the process converges in law to the diffusion process (our $X_t$)  defined by $ a $, $ b $ and thus and $s$ the jump process of generator $ Q $, since they are both solutions of the same problem of martingale. The proof of uniqueness is quite standard, following  Ethier \cite{ETHIER} when $ s $, $ m $, $ p $ constant. A straightforward modification allows to obtain the result for $ s $, $ m $, $ p $ random.

\subsection{Proof of Proposition \ref{prop:absorb}}
\label{absorb}
\quad\\
Let $g$ be the solution of the  differential equation \eqref{g}. Let us first verify that $ g$ is well defined on $[0,1]$. It must be ensured that the solutions do not diverge in 0 and 1, in which case the second member of the equation is not defined. For that we can write the solution of this equation. So $$g(x)=\int_{0}^{x}e^{-su}\left(K+\int_{\frac{1}{2}}^{u}\frac{e^{st}}{t(1-t)}dt\right)du+C$$ where $C,K$are constant. As $e^{-su}$ is bounded on$[0,1]$, there are two positive constants $B$ and $D$ such that $$\lim\limits_{x\rightarrow 0}|g(x)|\leqslant \lim\limits_{x\rightarrow 0}\int_{a}^{x}\int_{a}^{u} \frac{D}{t}dtdu + B\leqslant \infty$$
 as $\ln$ is integrable on a neighbourhood of 0. Thus $g$ is well defined on  [0,1] and bounded (because continue).\\

So we have  $g(X_{t\wedge T_{1, 0}})=g(X_{0})+\int_{0}^{t\wedge T_{1, 0}}g'(X_{u})\sqrt[]{2X_{u}(1-X_{u})}dB_{u}-t\wedge T_{1, 0}$.\\
But the process $\int_{0}^{t\wedge T_{1, 0}}g'(X_{u})\sqrt[]{2X_{u}(1-X_{u})}dB_{u}$ is a stopped martingale  because $t\wedge T_{1, 0}$
is a stopping time and $g'(X_{u})\sqrt[]{2X_{u}(1-X_{u})}$ is adapted to the considered filtration.
We deduce that $\E_{X_{0}}[g(X_{t\wedge T_{1, 0}})]=g(X_{0})-\E_{X_{0}}[t\wedge T_{1, 0}]$ and the first property, 
 i.e $\E_{X_{0}}[t\wedge T_{1, 0}]\leqslant 2 \sup\limits_{[0,1]}(g)\leqslant \infty$, and thus the second point is shown.Now if  $t\rightarrow \infty$, $\E_{X_{0}}[g(X_{t\wedge T_{1, 0}})]\rightarrow0$  because  $g(0)=g(1)=0$ and we find again $g(X_{0})=E_{X_{0}}[ T_{1, 0}]$.\\
To prove the third point, consider $f(x)=e^{-sx}-1$. Then $f$ is solution of 
$f''(x)+sf'(x)=0$ and  
$f(0)=0$. By It\^o's formula, we obtain $d(f(X_{t})=f(X_{0})+dMt $. As  $T_{1}$ and  $T_{0}$ are stopping times $T_{1}\wedge T_{0}=T_{1, 0}$ is also a stopping time. So$ f(X_{ T_{0,1} })=f(X_{0})+dM'_{t}$ où $dM'_{t}$ is still a martingale. By taking expectation we have $\E_{X_{0}}[f(X_{ T_{0,1}})]=f(X_{0})=f(1)\mathbb{P}(T_{1}<T_{0})$  and we deduce $ \mathbb{P}(T_{1}<T_{0})=\frac{e^{-sX_{0}}-1}{e^{-s}-1}$.

\subsection{Proof of Proposition \ref{prop:access}}
\label{access}
\quad\\
Let us consider the speed measure and the scale function as in Feller \cite{Feller2}.

\begin{align*}
m(y)&=\frac{1}{2y(1-y)}exp(\int_{a}^{y}s+\frac{m(p-x)}{x(1-x)}dx) \quad \textnormal{(speed measure)}\\ 
&=c\times y^{mp-1}\times (1-y)^{m(1-p)-1}\times \exp(sy),\quad c \in \mathbb{R}
\end{align*}
\begin{align*}
\mu(t)&=\int_{a}^{t} exp(-\int_{a}^{y}s+\frac{m(p-x)}{x(1-x)}dx)dy \quad \textnormal{(scale function)}\\
&=c'\int_{a}^{t} y^{-mp}\times (1-y)^{-m(1-p)}\times \exp(-sy)dy,\quad c' \in \mathbb{R}.
\end{align*}

Then $1$ is reachable if and only if  $\mu(1)<\infty$ and $\int_{\frac{1}{2}}^{1}\mu(1)-\mu(y)m(y)dy<\infty$. It is easily seen that $\mu(1)<\infty$ if and only if  $m(1-p)<1$.\\
Next,  $\int_{\frac{1}{2}}^{1}(\mu(1)-\mu(y)m(y))dy<\infty$  if and only if  $\mu(1)-\int\mu(y)m(y)$ is integrable on a neighborhood of 1. But 
\begin{align*}
\mu(1)-\mu(y)m(y)&\leqslant \int_{y}^{1}x^{-mp}(1-x)^{-m(1-p)}y^{mp-1}(1-y)^{m(1-p)-1}e^{sy-sx}dx\\
&\leqslant \int_{y}^{1}(1-x)^{-m(1-p)}dx\times(1-y)^{m(1-p)-1}\frac{1}{y}\\
&\leqslant a\times \frac{1}{y},
\end{align*}
for some constant $a$. This quantity is well defined and integrable on a neighborhood of  1. So 1 is reachable if and only if  $m(1-p)<1$.\\

Now if  $m(1-p) \leqslant 1$, $1$ is not reachable, it is regular (reflective barriers) if and only if $m(y)$ is  integrable. It is indeed the case here, $-21<m(1-p)-1\leqslant 0$. Of course the same holds for  $0$. 

\bibliographystyle{plain}
\bibliography{biblio}

\end{document}